\documentclass[a4paper,12pt]{article}
\usepackage{amsmath}
\usepackage{amssymb}
\usepackage{tabularx}
\usepackage{enumerate}
\usepackage{graphicx}
\usepackage{subfigure}
\usepackage{color}
\usepackage[pagewise]{lineno}
\usepackage{longtable}
\usepackage{epstopdf}

\usepackage[justification=centering]{caption}
\usepackage{multirow}
\usepackage{cases}
\usepackage[compress]{cite}
\usepackage{lineno}
\usepackage{color}
\usepackage{setspace}
\usepackage{multirow}
 \usepackage{rotating}
\usepackage{geometry}
\newtheorem{thm}{Theorem}[section]
\newtheorem{lem}[thm]{Lemma}

\newtheorem{exm}{Example}[section]
\newtheorem{prop}[thm]{Proposition}
\newtheorem{rmk}{Remark}[section]

\newtheorem{proof}{Proof}

\newcommand{\qed}{\hspace{1em}\mbox{\raisebox{0.65ex}{\fbox{}}}}

\numberwithin{equation}{section}

\newcommand{\be}{\begin{equation}}
\newcommand{\ee}{\end{equation}}
\newcommand\bes{\begin{eqnarray}} \newcommand\ees{\end{eqnarray}}
\newcommand{\bess}{\begin{eqnarray*}}
\newcommand{\eess}{\end{eqnarray*}}

\newcommand{\bpf}{{\bf Proof:\ \ }}
\newcommand{\epf}{\mbox{}\hfill $\Box$}

\begin{document}

\thispagestyle{empty}
\title{Spatial diffusion and periodic evolving of domain in an SIS epidemic model
\thanks{The work is partially supported by the NNSF of China (Grant No. 11771381).}}
\date{\empty}

\author{Yachun Tong$^{a}$ and Zhigui Lin$^a\thanks{Corresponding author. Email: zglin@yzu.edu.cn (Z. Lin).}$\\
{\small $^a$ School of Mathematical Science, Yangzhou University, Yangzhou 225002, China}\\
}

\maketitle
\begin{quote}
\noindent
{\bf Abstract.} { 
 \small In order to explore the impact of periodically evolving domain on the transmission of disease, we study
 a SIS reaction-diffusion model with logistic term on a periodically evolving domain. The basic reproduction
 number ${\mathcal{R}}_0$ is given by the next generation infection operator, and relies on the evolving rate
 of the periodically evolving domain, diffusion coefficient of infected individuals $d_I$ and size of the
 space. The monotonicity of ${\mathcal{R}}_0$ with respect to $d_I$, evolving rate $\rho(t)$ and interval length
 $L$ are derived, and asymptotic property of ${\mathcal{R}}_0$ if $d_I$ or $L$ is small enough or large enough in one-dimensional
 space are discussed. ${\mathcal{R}}_0$ as threshold can be used to characterize whether the disease-free equilibrium is
 stable or not. Our theoretical results and numerical simulations indicate that small evolving rate, small
 diffusion of infected individuals and small interval length have positive impact on prevention and control
 of disease .
}


\end{quote}

\section{Introduction}

In recent decades, it is generally recognized that spatial diffusion and environmental heterogeneity have
an important impact on the spread of many diseases. Therefore, more and more work is devoted to studying the
impact of diffusion on disease control and transmission \cite{AL, BB, cl, DW}.
In order to capture the impact of spatial heterogeneity of environment and movement of individuals on the
persistence and extinction of a disease, Allen et al. \cite {AL} proposed a frequency-dependent SIS
(susceptible-infected-susceptible) epidemic reaction-diffusion system
\begin{eqnarray}
\left\{
\begin{array}{lll}
S_{t}-d_S\Delta S=-\frac{\beta (x) SI}{S+I}+\gamma(x)I,\; &\, x\in\Omega,\ t>0, \\[2mm]
I_{t}-d_I\Delta I=\frac{\beta (x) SI}{S+I}-\gamma(x)I,\; &\, x\in\Omega,\ t>0,
\end{array} \right.
\label{a01}
\end{eqnarray}
where $S(x,t)$ and $I(x,t)$  represent the numbers of susceptible and infected individuals at location
$x$ and time $t$, respectively; $\beta(x)$ and $\gamma(x)$ are positive bounded H$\ddot{o}$lder continuous
functions on $\overline{\Omega}$, and denote the rates of contact transmission and disease recovery at $x$,
respectively; $d_S$ and $d_I$ are positive constants that denote the diffusion coefficient of susceptible
and infected individuals, respectively. Based on the basic reproduction number ${\mathcal{R}}_0$, the authors
characterised the risk of region. For high domain (${\mathcal{R}}_0>1$), the DFE (disease-free equilibrium) is
always unstable and there exists a unique EE (endemic equilibrium); in the low risk domain (${\mathcal{R}}_0<1$),
the DFE is stable if and only if infected individuals have mobility above a threshold value. For some special cases,
Peng and Liu \cite{pli} proved that unique EE is globally stable, which was conjecture proposed by
Allen, et al in \cite{AL}. Further results that the effect of individual movement (large or small) on the
existence and disappearance of disease were obtained in \cite{p}. Moreover, \cite{pz} discussed problem
\eqref{a01} with $\beta$ and $\gamma$ being functions of spatiotemporal variables and temporally periodic,
and discovered that spatial heterogeneity and temporal periodicity can enhance the persistence of the disease.
More work on system \eqref{a01} can be found in the literature \cite{dhk,cl,ctz,dhl,wz}.

Recently, Li, et al \cite{lht} considered the following SIS epidemic with logistic source
\begin{equation*}
\left\{ \begin{array}{llll}
S_{t}-d_S\Delta S=a(x)S-b(x)S^2-\beta(x)f(S,I)I+\gamma(x)I,&x\in\Omega,\,t>0,\\
I_{t}-d_I\Delta I=\beta(x)f(S,I)I-\gamma(x)I,&x\in\Omega,\,t>0,\\
\frac{\partial S}{\partial \nu}=\frac{\partial I}{\partial \nu}=0,&x\in\partial\Omega,\,t>0,\\
S(x,0)=S_0(x),\ I(x,0)= I_0(x),&x\in\Omega,
\end{array}\right.
\end{equation*}
where $f(S,I)= \frac{S}{S+I}$, $\beta(x)$, $\gamma(x)$, $a(x)$ and $b(x)$ are positive H\"{o}lder continuous
functions on $\overline \Omega$. The nonlinear term $a(x)S-b(x)S^2$ denotes that the susceptible population
is subject to logistic growth. They studied the stability of DFE and EE, and presented the asymptotic profiles
of EE as the motility of susceptible or infected individuals is small or large.

It is not difficult to find that all of the above SIS epidemic models are considered in fixed domain, that is,
the domain is independent of time $t$. In fact, in the natural world, biological habitats usually change over the time.
Mathematically, the free boundary can be used to describe the domain change induced by the movement of species.
Therefore, the free boundary problem is widely studied in many fields, for example, wound healing \cite{cf},
the spread of invasive species \cite{DL, GW, WZ1, WZ2}, the transmission of infected disease \cite{CLW, CLY, GKLZ, KLZ},
and so on. However, the free boundary is unknown. On the other hand, the evolving domain is known, such as the change
of lake depth due to the alternation of seasons \cite{WKL}, and the spread of mosquito borne diseases due to
climate \cite{bam}. Owing to the turnover of the four seasons, the evolving of domain is usually periodic. In recent years, some scholars
have studied mathematical models on a periodically evolving domain, for example, an SIS epidemic model \cite{pl},
a mutualistic model \cite{alt}, an Aedes aegypti mosquito model \cite{zl}, a dengue fever model \cite{zxc}, and
so on. In this paper, we consider the following problem
\begin{equation}
\left\{ \begin{array}{llll}
S_{t}-d_S\Delta S=a(x)S-b(x)S^2-\beta(x)f(S,I)I+\gamma(x)I,&x\in\Omega,\,t>0,\\
I_{t}-d_I\Delta I=\beta(x)f(S,I)I-\gamma(x)I,&x\in\Omega,\,t>0,\\
\frac{\partial S}{\partial \nu}=\frac{\partial I}{\partial \nu}=0,&x\in\partial\Omega,\,t>0,\\
S(x,0)=S_0(x),\ I(x,0)= I_0(x),&x\in\Omega
\end{array}\right.
\label{a02}
\end{equation}
on a periodically evolving domain.

As in \cite{C1, tql1}, let $\Omega (t)\subset \mathbf R^n(n\geq1)$ be a simply connected bounded shifting
domain at time $t\ge 0$ with its evolving boundary $\partial\Omega(t)$, and then $\bar{S}(x(t),t)$
and $\bar{I}(x(t),t)$ are the densities of susceptible and infected species at position $x(t)$ and
time $t\ge 0$. We can establish the evolution equation on a evolving domain
by conservation of mass and Reynolds transport theorem \cite{tql2}. For any point
$$
x(t)=(x_1(t),x_2(t),  \ldots, x_n(t))  \in \Omega(t).
$$
The evolution of the domain produces a velocity of flow $\mathbf a=\dot{x}(t)$. By Reynolds
transport theorem \cite{tql2}, we have
\begin{eqnarray}
\begin{array}{lll}
\frac{\partial \bar{S}}{\partial t}+\nabla \bar{S}\cdot \mathbf a+\bar{S}(\nabla\cdot \mathbf a)=
d_S \Delta \bar{S}+ f_1(\bar{S},\bar{I},t)  \ \; \textrm  {in}\  \Omega(t), \\[2mm]
\frac{\partial \bar{I}}{\partial t}+\nabla \bar{I}\cdot \mathbf a+\bar{I}(\nabla\cdot \mathbf a)=
d_I \Delta \bar{I}+ f_2(\bar{S},\bar{I},t)  \ \; \textrm  {in}\  \Omega(t), \\[2mm]
\end{array}
\label{a03}
\end{eqnarray}
where
$$f_1(\bar{S},\bar{I},t)=a(x)\bar{S}-b(x)\bar{S}^2-\beta(x)f(\bar{S},\bar{I})\bar{I}+\gamma(x)\bar{I},\;
f_2(\bar{S},\bar{I},t)=\beta(x)f(\bar{S},\bar{I})\bar{I}-\gamma(x)\bar{I},$$
$\nabla \bar{S} \cdot \mathbf a$
and $\nabla \bar{I} \cdot \mathbf a$ are called advection terms generated by the evolving domain while
$(\nabla \cdot\mathbf a)\bar{S} $ and $(\nabla \cdot\mathbf a)\bar{I} $ are called dilution terms caused
by the change of element volume and local volume due to local growth and flow movement \cite{CGM, CGM2}.
In order to deal with the convection term and dilution term caused by regional evolution in equation
\eqref{a03}, we use Lagrangian transformations to change equations \eqref{a03}. Let $y_1, y_2, \ldots, y_n$
be fixed cartesian coordinates in a fixed domain $\Omega(0)$ such that
$$
x_1(t)=\hat x_1(y_1, y_2, \ldots, y_n, t),
$$
$$
x_2(t)=\hat x_2(y_1, y_2, \ldots, y_n, t),
$$
$$
\dots
$$
$$
x_n(t)=\hat x_n(y_1, y_2, \ldots, y_n, t).
$$
In this transformation, $(\bar{S}, \bar{I})$ is then mapped into the new vector $(S, I)$ defined as
\begin{eqnarray}
\begin{array}{lll}
\bar{S}(x_1(t), x_2(t), \ldots, x_n(t), t)=S(y_1,y_2,\ldots, y_n, t),\\[2mm]
\bar{I}(x_1(t), x_2(t), \ldots, x_n(t), t)=I(y_1,y_2,\ldots, y_n, t).
\end{array}
\label{a04}
\end{eqnarray}
Now equations (\ref{a03}) is translated into another form which are defined on the fixed domain $\Omega(0)$
with respect to $y=(y_1,y_2,\ldots, y_n)$, but the new equations are very complicated for arbitrary domain
evolution. To further simplify equations (\ref{a03}), we all suppose that domain evolution is uniform
and isotropic. In other words, the evolution of the domain occurs in the same proportion in all directions
as time increases. From a mathematical point of view, $x(t)=(x_1(t),x_2(t),\ldots, x_n(t))$ can be described
as:
\bes
(x_1(t),x_2(t),\ldots, x_n(t))=\rho(t)(y_1,y_2,\ldots, y_n),
&y\in \Omega(0),
\label{a05}
\ees
where the positive continuous function $\rho(t)$ is called evolving rate subject to $\rho(0)=1$. The domain is
periodically evolving if $\rho(t)=\rho(t+T)$ for some $T>0$ \cite{wzc, SSM}. The domain is called growing one
If $\dot \rho(t) \ge0$ \cite{tql2, tql1}, and when $\dot \rho(t) \le0$, the domain is shrinking \cite{WN}.
Using \eqref{a05} yields
\bess
& & S_t=\bar{S}_t+\nabla \bar{S} \cdot \mathbf a,\quad I_t=\bar{I}_t+\nabla \bar{I} \cdot \mathbf a,\\[2mm]
& & \mathbf a=\dot{x}(t)=\dot\rho(t)(y_1, y_2,\ldots, y_n)=\frac{\dot\rho(t)}{\rho(t)}(x_1, x_2, \ldots, x_n), \\[2mm]
& & \nabla \cdot \mathbf a=\frac{n \dot\rho(t)}{\rho(t)},\quad
\Delta \bar{S}=\frac{1}{\rho^2(t)}\Delta S,\quad
\Delta \bar{I}=\frac{1}{\rho^2(t)}\Delta I.
\eess
Then \eqref{a03} becomes
\begin{eqnarray}
\left\{
\begin{array}{lll}
S_t=\displaystyle\frac{d_S}{\rho^2(t)}\Delta S-\frac{n\dot \rho(t)}{\rho(t)}S + f_1(S,I,t), & y \in \Omega(0),\,\ t>0,\\[3mm]
I_t=\displaystyle\frac{d_I}{\rho^2(t)}\Delta I-\frac{n\dot \rho(t)}{\rho(t)}I + f_2(S,I,t), & y \in \Omega(0),\,\ t>0.
\end{array} \right.
\label{a06}
\end{eqnarray}

Now we transform the SIS epidemic model on a periodically evolving domain $\Omega(t)$ into the
following problem on a fixed domain $\Omega(0)$:
\begin{small}
\begin{eqnarray}
\left\{
\begin{array}{llllll}
S_{t}-\frac{d_S}{\rho^2(t)} \Delta S=a(\rho(t)y)S-b(\rho(t)y)S^{2}-\beta(\rho(t)y)f(S,I)I & \\
\qquad \qquad \qquad +\gamma(\rho(t)y)I-
\frac{n\dot \rho(t)}{\rho(t)}S,\; & y\in \Omega(0), \; t>0,\\[2mm]
I_{t}-\frac{d_I}{\rho^2(t)} \Delta I=\beta(\rho(t)y)f(S,I)I-\gamma(\rho(t)y)I-\frac{n\dot \rho(t)}{\rho(t)}I,
\;& y\in \Omega(0),\;  t>0, \\[2mm]
\frac{\partial S}{\partial \nu}=\frac{\partial I}{\partial \nu}=0,  &y\in \partial \Omega(0),\;  t>0,
\end{array} \right.
\label{a07}
\end{eqnarray}
\end{small}
with the initial condition
\begin{equation}
S(y,0)=S_0(y)>0,\, I(y,0)=I_0(y)\geq0,\ I_0(y)\not\equiv 0,\  y\in  \overline{\Omega(0)}.
\label{a08}
\end{equation}
For later application, we also consider problem \eqref{a07} with the periodic condition
\begin{equation}
S(y,0)=S(y,T),\, I(y,0)=I(y,T),\; \ y\in \overline{\Omega(0)},
\label{a09}
\end{equation}
where $f(S,I)=\frac {S}{S+I}$ is monotonically decreasing with respect to $I$ and increasing with
respect to $S$ and $\lim\limits_{I\rightarrow0}f(S,I)=1$.

The remainder of the article is organized as below. Section 2 is devoted to the existence and uniqueness
of DFE, and Section 3 gives the definition of basic reproduction number and its influence
on the stability of DFE. Section 4 is concerned with the limiting behaviors of ${\mathcal{R}}_0$
with respect to $d_I$ and $L$, and monotonicity of ${\mathcal{R}}_0$ with respect to $d_I$, $L$
and $\rho(t)$. In the last section, we explain our theoretical results by the numerical simulations.

\section{\bf The existence and uniqueness of DFE}

In this section, We focus on the existence and uniqueness of the disease-free equilibrium $(S^{*}(y,t),0)$,
where $(S^{*}(y,t),0)$ satisfies
\begin{eqnarray}
\left\{
\begin{array}{lll}
S_{t}-\frac{d_S}{\rho^2(t)} \Delta S=a(\rho(t)y)S-b(\rho(t)y)S^{2}-\frac{n\dot \rho(t)}{\rho(t)}S,\; &\ y\in
\Omega(0), 0<t\leq T,\\[2mm]
\frac{\partial S}{\partial \nu}=0, \; &\ y\in \partial\Omega(0), 0<t<T,\\[2mm]
S(y,0)=S(y,T),\; &\ y\in \overline{\Omega(0)}.
\end{array} \right.
\label{b01}
\end{eqnarray}

For existence of the nontrivial solution, we consider the following T-periodic parabolic eigenvalue problem
\begin{eqnarray}
\left\{
\begin{array}{lll}
\psi_{t}-\frac{d_S}{\rho^2(t)} \Delta \psi=a(\rho(t)y)\psi-\frac{n\dot \rho(t)}{\rho(t)}\psi+\lambda\psi,\; &\
 y\in \Omega(0), 0<t\leq T,\\[2mm]
\frac{\partial \psi(y,t)}{\partial \nu}=0,\; &\ y\in \partial\Omega(0), 0<t<T,\\[2mm]
\psi(y,0)=\psi(y,T),\; &\ y\in \overline{\Omega(0)}.
\end{array} \right.
\label{b02}
\end{eqnarray}
Let $\lambda(a)$ is the principal eigenvalue of the above the eigenvalue problem \eqref{b02} and $\psi^{*}>0$ to
be an eigenfunction corresponding to $\lambda(a)$. Now if we multiply \eqref{b02} by $\frac{1}{\psi^{*}}$ and
then integrate the first equation of \eqref{b02} over $\Omega(0)\times[0,T]$, we are led to
$$\begin{array}{llllll}
&&\int_{0}^{T}\int_{\Omega(0)}\frac{\psi^{*}_{t}}{\psi^{*}}dtdy-\int_{0}^{T}\int_{\Omega(0)}\frac{d_S}{\rho^2(t)}
\frac{\Delta \psi^{*}}{\psi^{*}}dtdy\\[2mm]
&=&\int_{0}^{T}\int_{\Omega(0)}a(\rho(t)y)dtdy-\int_{0}^{T}\int_{\Omega(0)}\frac{n\dot \rho(t)}{\rho(t)}dtdy+
\int_{0}^{T}\int_{\Omega(0)}\lambda(a)dtdy.
\end{array}$$
Recalling that $a(\rho(t)y)>0$ and $d_{S}>0$, we have
\begin{eqnarray}
\lambda(a)=\frac{-\int_{0}^{T}\int_{\Omega(0)}\frac{d_S}{\rho^2(t)} \frac{|\nabla \psi^{*}|^{2}}{{\psi^{*}}^{2}}dtdy-\int_{0}^{T}\int_{\Omega(0)}a(\rho(t)y)dtdy}{|\Omega(0)|\times T}<0.
\label{b03}
\end{eqnarray}

Taking $\underline{S}=\varepsilon\psi^{*}$, then
$$\begin{array}{llll}
&&\underline{S}_{t}-\frac{d_S}{\rho^2(t)} \Delta\underline{S}-a(\rho(t)y)\underline{S}+
b(\rho(t)y)\underline{S}^{2}+\frac{n\dot \rho(t)}{\rho(t)}\underline{S}\\[2mm]
&=&\varepsilon a(\rho(t)y)\psi^{*}-\varepsilon\frac{n\dot\rho(t)}{\rho(t)}\psi^{*}+
\lambda(a)\varepsilon\psi^{*}-\varepsilon a(\rho(t)y)\psi^{*}+b(\rho(t)y)\varepsilon^{2}{\psi^{*}}^{2}
+\varepsilon\frac{n\dot \rho(t)}{\rho(t)}\psi^{*}\\[2mm]
&=&\lambda(a)\varepsilon\psi^{*}+b(\rho(t)y)\varepsilon^{2}{\psi^{*}}^{2}\\[2mm]
&=&\varepsilon\psi^{*}(\lambda(a)+b(\rho(t)y)\varepsilon\psi^{*}).
\end{array}$$
Due to \eqref{b03}, $\underline{S}$ is the lower solution of problem \eqref{b01} if $\varepsilon$ is small enough
such that $\lambda(a)+b(\rho(t)y)\varepsilon\psi^{*}<0$. Obviously, $\overline{S}=M$ is the upper solution
of problem \eqref{b01} when positive constant $M$ is sufficiently large. According to the upper and lower solutions
method, problem \eqref{b01} has at least one positive solution $S^*(y,t)$ satisfying $\underline{S}\leq S^*\leq \overline{S}$
for $(y,t)\in\overline{\Omega(0)}\times [0,T]$.

Finally, we prove the uniqueness of the solution to problem \eqref{b01}. Considering the following auxiliary problem
\begin{eqnarray}
\left\{
\begin{array}{lll}
u_{t}-\frac{d_S}{\rho^2(t)} \Delta u=a(\rho(t)y)u-b(\rho(t)y)u^{2}-\frac{n\dot \rho(t)}{\rho(t)}u,\; &\
 y\in \Omega(0), t>0,\\[2mm]
\frac{\partial u(y,t)}{\partial \nu}=0,\; &\ y\in \partial\Omega(0), t>0,\\[2mm]
u(y,0)=u_{0}(y),\; &\ y\in \overline{\Omega(0)},
\end{array} \right.
\label{b04}
\end{eqnarray}
where for any $u_{0}\in C^{1}(\overline{\Omega(0)})\cap W^{2} _{p}(\Omega(0))$ satisfying
$\frac{\partial u_{0}}{\partial \nu}=0$ on $\partial\Omega(0)$ and $u_{0}>0$ in $\Omega(0)$.
We can easily to know that the initial-boundary value problem \eqref{b04} has an unique positive solution.

Let $S_{1}$ and $S_{2}$ be two positive solutions of problem \eqref{b01}. Invoking the upper and lower
solutions method, there exists a $M\gg1$ such that $M\geq S_{j}$, $j=1, 2$. Let $u$ be the unique solution
of problem \eqref{b04} with initial data $u_{0}(y)=M$. By the comparison principle, we have $M\geq u\geq S_{j}$.
Define
$$
u^{i}(y,t)=u(y,t+iT), \ \ \  \ (y,t)\in \Omega(0)\times[0,T).
$$
Due to that $a$ and $b$ are T-periodic in the time, so $u^{i}$ satisfies
\begin{equation*}
\left\{
\begin{array}{lll}
 u_{t}^{i}-\frac{d_S}{\rho^2(t)} \Delta u^{i}=a(\rho(t)y)u^{i}-b(\rho(t)y)(u^{i})^{2}-
 \frac{n\dot \rho(t)}{\rho(t)}u^{i},\; &\ y\in \Omega(0), 0<t\leq T,\\[2mm]
 \frac{\partial u^{i}(y,t)}{\partial \nu}=0,\; &\ y\in \partial\Omega(0), 0<t<T,\\[2mm]
 u^{i}(y,0)=u(y,iT),\; &\ y\in \overline{\Omega(0)}.
\end{array} \right.
\end{equation*}
Since
$$
u(y,0)=M\geq u(y,T)=u^{1}(y,0)\geq S_{j}(y,T)=S_{j}(y,0).
$$
We can use the comparison principle to derive $M\geq u\geq u^{1}\geq S_{j}$ on $\overline{\Omega(0)}\times[0,T]$,
which in turn implies
$$
u^{1}(y,0)=u(y,T)\geq u^{1}(y,T)=u^{2}(y,0)\geq S_{j}(y,T)=S_{j}(y,0).
$$
By the above inequality, $u^{1}\geq u^{2}\geq S_{j}$ on $\overline{\Omega(0)}\times[0,T]$.
Repeat the last step, we can obtain $u^{i}$ is monotone decreasing in $i$ and $M\geq u^{i}\geq S_{j}$ on
$\overline{\Omega(0)}\times[0,T]$. Hence, there exists a positive function $S_{M}$ such that
$$
u^{i}\rightarrow S_{M} \ \mbox{pointwise on}\ \overline{\Omega(0)}\times[0,T] \ \ \mbox{as} \ \ i\rightarrow\infty.
$$
Since $u^{i+1}(y,0)=u^{i}(y,T)$, so $S_{M}(y,0)=S_{M}(y,T)$. According to the regularity theory
and compactness argument, we know that
$$
u^{i}\rightarrow S_{M} \ \mbox{in}\ C^{2,1}(\overline{\Omega(0)}\times[0,T])
$$
and $M\geq S_{M} \geq S_{j}$ on $\overline{\Omega(0)}\times[0,T]$. Moreover,
$S_{M}$ is positive solution of problem \eqref{b01}.
Next, we show $S_{M}=S_{j},j=1,2$, which implies problem \eqref{b01} has an unique positive solution.
Without loss of generality, we just prove $S_{M}=S_{1}$. Since $S_{1}\leq S_{M}$, there
exists a positive constant $\theta$ such that $S_{1}\geq\theta S_{M}$ on $\overline{\Omega(0)}\times[0,T]$.
Setting
$$
\mu=\sup\{\theta: S_1\geq\theta S_M \  \textrm{in} \ \ \overline{\Omega(0)}\times[0,T]\}.
$$
obviously, $\mu\leq1$ and $S_{1}(y,t)\geq\mu S_{M}(y,t)$ on $\overline{\Omega(0)}\times[0,T]$.
If we can prove $\mu=1$, then $S_{1}=S_{M}$. We prove it by contradiction. Assume $\mu<1$, define
$\varphi=S_{1}-\mu S_{M}$. Then $\varphi\geq0$ and
$$
\frac{\partial \varphi}{\partial \nu}=0 \ \ \mbox{in} \ \ \partial\Omega(0)\times[0,T], \ \varphi(y,0)=
\varphi(y,T) \ \mbox{on} \ \overline{\Omega(0)}\times[0,T].
$$
Straightforward calculation can be obtained
$$\begin{array}{lll}
&&\varphi_{t}-\frac{d_S}{\rho^2(t)} \Delta\varphi=a(\rho(t)y)\varphi-b(\rho(t)y)({S_{1}}^{2}-\mu{S_{M}}^{2})
-\frac{n\dot \rho(t)}{\rho(t)}\varphi\\[2mm]
&>&a(\rho(t)y)\varphi-b(\rho(t)y)({S_{1}}^{2}-\mu^{2}{S_{M}}^{2}) -\frac{n\dot \rho(t)}{\rho(t)}\varphi\\[2mm]
&\geq&a(\rho(t)y)\varphi-2b(\rho(t)y)S_{1}\varphi-\frac{n\dot \rho(t)}{\rho(t)}\varphi.
\end{array}$$
By the maximum principle, $\varphi>0$ in $\Omega(0)\times[0,T]$. Hence, there exists $\tau>0$ such
that $\varphi\geq\tau S_{M}$. Therefore, $S_{1}\geq(\mu+\tau)S_{M}$ on $\overline{\Omega(0)}\times[0,T]$.
This is a contradiction with the definition of $\mu$. So $S_{M}=S_{1}$.

Hence, problem \eqref{b01} admits an unique positive periodic solution.

\section{\bf Threshold dynamics in terms of ${\mathcal{R}}_0$}
In this section, we first define the basic reproduction number ${\mathcal{R}}_0$, and then we study threshold
dynamics in terms of ${\mathcal{R}}_0$ for problem \eqref{a07}-\eqref{a08}. In general, the basic reproductive
number is used to describe the transmission mechanism of the diseases. Biologically,
the basic reproduction number ${\mathcal{R}}_0$ is calculated by the next generation matrix method \cite{DW} for epidemic models described by spatially-independent systems, or
calculated by the next infection operator \cite{WZ} for models described by spatially-dependent systems.

First, linearizing problem \eqref{a07} at the disease-free equilibrium $(S^*(y,t),0)$ gives
\begin{equation*}
\left\{
\begin{array}{ll}
U_t-D(t)\Delta U=A(t)U -B(t)U,\; &y\in \Omega(0), t>0, \\[2mm]
\frac{\partial F(y,t)}{\partial \nu}=0, &y\in \partial \Omega(0), t>0,
\end{array} \right.
\end{equation*}
with the same periodic condition \eqref{a09}, where
$$
U=\left(
    \begin{array}{c}
      u \\
      v \\
    \end{array}
  \right), \
  D(t)=\left(
         \begin{array}{cc}
           \frac{d_S}{\rho^2(t)} & 0 \\
           0 & \frac{d_I}{\rho^2(t)} \\
         \end{array}
       \right), \
       A(t)=\left(
       \begin{array}{cc}
         a(\rho(t)y) & \gamma(\rho(t)y) \\
         0 & \beta(\rho(t)y) \\
       \end{array}
     \right),
$$
$$
B(t)=\left(
       \begin{array}{cc}
         \frac{n\dot{\rho}(t)}{\rho(t)}+2b(\rho(t)y)S^{*}(\rho(t)y) & \beta(\rho(t)y) \\
         0 & \gamma(\rho(t)y)+\frac{n\dot{\rho}(t)}{\rho(t)} \\
       \end{array}
     \right).
$$
Let $V(t, s)$ be the evolution operator of the problem
\begin{equation*}
\left\{
\begin{array}{ll}
w_t-D(t)\Delta w=-B(t)w,\; &y\in \Omega(0), \;  t>0, \\[2mm]
\frac{\partial w(y,t)}{\partial \nu}=0, &y\in \partial \Omega(0),  \; t>0.
\end{array} \right.
\end{equation*}
By the standard theory of evolution operators, there exist positive constants $K$ and $c_0$ such that
$$\|V(t,s)\|\leq Ke^{-c_0(t-s)},\ \ \  \forall t\geq s,\, t, s\in R.$$
The ordered Banach space $C_T$ is made up of all $T-$ periodic and continuous function from $R$ to
$C(\overline{\Omega}(0),R)$, which is equipped with the maximum norm $\|\cdot\|$ and the positive cone
$C_T^+:=\{\xi\in C_T:\xi(t)y\geq0, \forall \ t\in R, y\in\overline{\Omega}(0)\}$. Set $\xi(y,t):=\xi(t)y$
for any given $\xi\in C_T$. Assume that $\eta=(\xi, \zeta)\in C_T\times C_T$ represents the density
distribution of $w$ at the spatial locaton $y\in \Omega(0)$ and time $s$. We introduce the linear operator
$$
L(\eta)(t):=\int_{0}^{\infty}V(t,t-s)A(\cdot,t-s)\eta(\cdot,t-s)ds,
$$
which is called as the next infection operator in \cite{wwz}.
It is clear that $L$ is positive, continuous and compact on $C_T\times C_T$ under our assumption on
$a(\rho(t)y)$, $b(\rho(t)y)$, $\beta(\rho(t)y)$, $\gamma(\rho(t)y)$ and $S^{*}(\rho(t)y)$. We define
the spectral radius of $L$
$$
{\mathcal{R}}_0=r(L)
$$
as the basic reproduction number for periodic system \eqref{a07}, \eqref{a09}. Furthermore, we have
the following results.

\begin{lem}
$(i)$  Let $\mu_0$ is the principal eigenvalue of the following periodic-parabolic eigenvalue problem
\begin{eqnarray}
\left\{
\begin{array}{lll}
\Psi_{t}-\frac{d_S}{\rho^2(t)} \Delta\Psi=\frac{\gamma(\rho(t)y)}{\mu}\Phi+a(\rho(t)y)\Psi-2b(\rho(t)y)S^{*}\Psi & \\
\qquad \qquad \qquad \quad -\beta(\rho(t)y)\Phi-\frac{n\dot \rho(t)}{\rho(t)}\Psi,\; &\ y\in \Omega(0), 0<t\leq T, \\[2mm]
\Phi_{t}-\frac{d_I}{\rho^2(t)} \Delta\Phi=\frac{\beta(\rho(t)y)}{\mu}\Phi-\gamma(\rho(t)y)\Phi-
\frac{n\dot \rho(t)}{\rho(t)}\Phi,\; &\ y\in \Omega(0), 0<t\leq T, \\[2mm]
\frac{\partial \Psi(y,t)}{\partial \nu}=\frac{\partial \Phi(y,t)}{\partial \nu}=0,\ & \
y\in \partial \Omega(0), 0<t\leq T,\\[2mm]
\Psi(y,0)=\Psi(y,T),\, \Phi(y,0)=\Phi(y,T),\; &\ y\in \overline{\Omega(0)},
\end{array} \right.
\label{c01}
\end{eqnarray}
then, ${\mathcal{R}}_0=\mu_0$.

$(ii)\  \rm{sign} (1-{\mathcal{R}}_0)=\rm{sign}\lambda_0$, where $\lambda_0$ is the principal eigenvalue of
the following reaction-diffusion problem
\begin{small}
\begin{eqnarray}
\left\{
\begin{array}{lll}
\Psi_{t}-\frac{d_S}{\rho^2(t)} \Delta\Psi=[\gamma(\rho(t)y)-\beta(\rho(t)y)]\Phi+\lambda\Psi &\\
\qquad \qquad \qquad \quad +[a(\rho(t)y)-2b(\rho(t)y)S^*-\frac{n\dot \rho(t)}{\rho(t)}]\Psi, \; &\ y\in \Omega(0), 0<t\leq T, \\[2mm]
\Phi_{t}-\frac{d_I}{\rho^2(t)} \Delta\Phi=\beta(\rho(t)y)\Phi-\gamma(\rho(t)y)\Phi-\frac{n\dot \rho(t)}{\rho(t)}\Phi+
\lambda\Phi,\; &\ y\in \Omega(0), 0<t\leq T, \\[2mm]
\frac{\partial \Psi(y,t)}{\partial \nu}=\frac{\partial \Phi(y,t)}{\partial \nu}=0,\ & \ y\in \partial \Omega(0), 0<t\leq T,\\[2mm]
\Psi(y,0)=\Psi(y,T),\, \Phi(y,0)=\Phi(y,T),\; &\ y\in \overline{\Omega(0)}.
\end{array} \right.
\label{c02}
\end{eqnarray}
\end{small}
\end{lem}

From \eqref{c01}, $({\mathcal{R}}_0,\Phi)$ satisfies the following problem
\begin{eqnarray*}
\left\{
\begin{array}{lll}
\Phi_{t}-\frac{d_I}{\rho^2(t)} \Delta\Phi=\frac{\beta(\rho(t)y)}{{\mathcal{R}}_0}\Phi-\gamma(\rho(t)y)\Phi-
\frac{n\dot \rho(t)}{\rho(t)}\Phi,\; &\ y\in \Omega(0), 0<t\leq T, \\[2mm]
\frac{\partial \Phi(y,t)}{\partial \nu}=0,\ & \ y\in \partial \Omega(0), 0<t\leq T,\\[2mm]
 \Phi(y,0)=\Phi(y,T),\; &\ y\in \overline{\Omega(0)},
\end{array} \right.
\end{eqnarray*}
and $\Phi$ satisfies
\begin{small}
\begin{equation}
 \Phi_{t}-\frac{d_I}{\rho^2(t)} \Delta\Phi+[\max\limits_{y\in\overline{\Omega(0)}}\gamma(\rho(t)y)+\frac{n\dot \rho(t)}{\rho(t)}]\Phi\geq\frac{\beta(\rho(t)y)}{{\mathcal{R}}_0}\Phi\geq
 \frac{\min\limits_{y\in\overline{\Omega(0)}}\beta(\rho(t)y)}{{\mathcal{R}}_0}\Phi
\label{c03}
\end{equation}
\end{small}
for $(y,t)\in \overline{\Omega(0)}\times[0,T]$.
Now, we consider the following eigenvalue problem
{\small \begin{eqnarray}
\left\{
\begin{array}{ll}
W_{t}-\frac{d_I}{\rho^2(t)} \Delta W+[\max\limits_{y\in\overline{\Omega(0)}}\gamma(\rho(t)y)+\frac{n\dot \rho(t)}{\rho(t)}]W=\frac{\min\limits_{y\in\overline {\Omega(0)}}\beta(\rho(t)y)}{\mathcal{R}}W,\; &\
 y\in \Omega(0), 0<t\leq T, \\[2mm]
\frac{\partial W(y,t)}{\partial \nu}=0,\ & \ y\in \partial \Omega(0), 0<t\leq T,\\[2mm]
W(y,0)=W(y,T),\; &\ y\in \overline{\Omega(0)}.
\end{array} \right.
\label{c04}
\end{eqnarray}}
It is easy to find that the eigenvalue problem \eqref{c04} is equivalent to the ODE problem
\begin{eqnarray}
\left\{
\begin{array}{ll}
\frac{dW}{dt}+[\max\limits_{y\in\overline{\Omega(0)}}\gamma(\rho(t)y)+\frac{n\dot \rho(t)}{\rho(t)}]W=
\frac{\min\limits_{y\in\overline{\Omega(0)}}\beta(\rho(t)y)}{\mathcal{R}}W,\; &\ 0<t\leq T, \\[2mm]
W(0)=W(T).
\end{array} \right.
\label{c05}
\end{eqnarray}
We notice that \eqref{c05} has a positive solution $W(t)$ if and only if
$$
\mathcal{R}=\frac{\int_0^T\min\limits_{y\in\overline{\Omega(0)}}\beta(\rho(t)y)dt}{\int_0^T \max\limits_{y\in\overline{\Omega(0)}}\gamma(\rho(t)y)dt}.
$$
In light of Proposition 5.2 \cite{pz2}, \eqref{c03} and \eqref{c05}, we deduce that
\begin{equation}
{\mathcal{R}}_0\geq\mathcal{R}=\frac{\int_0^T\min\limits_{y\in\overline{\Omega(0)}}\beta(\rho(t)y)dt}
{\int_0^T \max\limits_{y\in\overline{\Omega(0)}}\gamma(\rho(t)y)dt}.
\label{c06}
\end{equation}
The similar analysis to above yields
\begin{equation}
{\mathcal{R}}_0\leq\frac{\int_0^T\max\limits_{y\in\overline{\Omega(0)}}\beta(\rho(t)y)dt}{\int_0^T\min\limits_{y\in\overline{\Omega(0)}}\gamma(\rho(t)y)dt}.
\label{c07}
\end{equation}

In light of \eqref{c06} and \eqref{c07}, one easily gain
\begin{equation}
\frac{\int_0^T\min\limits_{y\in\overline{\Omega(0)}}\beta(\rho(t)y)dt}{\int_0^T\max\limits_{y\in\overline {\Omega(0)}}\gamma(\rho(t)y)dt}\leq{\mathcal{R}}_0\leq\frac{\int_0^T\max\limits_{y\in\overline{\Omega(0)}}\beta(\rho(t)y)dt}{\int_0^T\min\limits_{y\in\overline {\Omega(0)}}\gamma(\rho(t)y)dt}.
\label{c08}
\end{equation}

Now we are going to study threshold dynamics in terms of ${\mathcal{R}}_0$ for problem \eqref{a07}-\eqref{a08},
which implies that ${\mathcal{R}}_0$ as threshold can be used to characterize stability of
the disease-free equilibrium.

\begin{thm}
\label{main}
The following statements hold:\\

$(i)$ If ${\mathcal{R}}_0<1$, then any nonnegative solutions of problem \eqref{a07}-\eqref{a08} satisfies
$\lim\limits_{t\to\infty}{I(y,t)}=0$ for $y\in\overline{\Omega(0)}$ and $\lim\limits_{m\rightarrow+\infty}S(y,t+mT)=
S^{*}(y,t)$ for $(y,t)\in \overline{\Omega(0)}\times[0,\infty)$. That is also to say, the disease-free equilibrium
$(S^{*}(y,t),0)$ is globally asymptotically stable;

$(ii)$ If ${\mathcal{R}}_0>1$, then there exists $\varepsilon_0>0$ such that any nonnegative nontrivial
solution of system \eqref{a07}-\eqref{a08} satisfies  $\limsup\limits_{t\to\infty}\|(S(y,t),I(y,t))-(S^*(y,t),0)\|\geq
\varepsilon_0$. That is to say, the disease-free equilibrium $(S^{*}(y,t),0)$ is unstable.
\end{thm}

\bpf  $(i)$ Suppose that $({\mathcal{R}}_0; \Psi, \Phi)$ is a solution of problem \eqref{c01}
and $\Psi(y,t), \Phi(y,t)>0$ for $(y,t)\in\overline{\Omega(0)}\times[0,T]$.
Setting $\bar{I}(y,t)=Me^{-\lambda t}\Phi(y,t)$, where $0<\lambda\leq \beta(\rho(t)y)(\frac{1}{{\mathcal{R}}_0}-1)$
for $(y,t)\in\overline{\Omega(0)}\times[0,T]$. According to $0\leq f(S,I)\leq 1$, then
$$\begin{array}{llllll}
&&\bar{I}_{t}-\frac{d_I}{\rho^2(t)} \Delta\bar{I}
-\beta(\rho(t)y)f(S,\bar{I})\bar{I}+\gamma(\rho(t)y)\bar{I}+\frac{n\dot \rho(t)}{\rho(t)}\bar{I}\\[2mm]
&\geq&\bar{I}_{t}-\frac{d_I}{\rho^2(t)} \Delta\bar{I}
-\beta(\rho(t)y)\bar{I}+\gamma(\rho(t)y)\bar{I}+\frac{n\dot \rho(t)}{\rho(t)}\bar{I}\\[2mm]
&=&Me^{-\lambda t}\Phi_t-\lambda Me^{-\lambda t}\Phi-Me^{-\lambda t}\frac{d_I}{\rho^2(t)} \Delta\Phi-\beta(\rho(t)y)
Me^{-\lambda t}\Phi\\[2mm]
&+&\gamma(\rho(t)y) Me^{-\lambda t}\Phi+ \frac{n\dot \rho(t)}{\rho(t)}Me^{-\lambda t}\Phi\\[2mm]
&=&\bar{I}\{-\lambda+\frac{\beta(\rho(t)y)}{{\mathcal{R}}_0}-\beta(\rho(t)y)\}\\[2mm]
&\geq& 0.
\end{array}$$
Hence, if $M$ is large enough, then $\bar{I}$ is the upper solution of the following problem
\begin{eqnarray*}
\left\{
\begin{array}{lllll}
{I}_{t}-\frac{d_I}{\rho^2(t)} \Delta{I}=\beta(\rho(t)y)f(S,{I}){I}-\gamma(\rho(t)y){I}-\frac{n\dot \rho(t)}{\rho(t)}{I},
\; &\ y\in \Omega(0), t>0, \\[2mm]
\frac{\partial I(y,t)}{\partial \nu}=0,\ & \ y\in \partial \Omega(0), t>0,\\[2mm]
I(y,0)=I_0(y)\geq0,\,I_0(y)\not\equiv0, \; &\ y\in  \overline{\Omega(0)}.
\end{array} \right.
\end{eqnarray*}
Since $\lim\limits_{t\rightarrow+\infty}\bar{I}(y,t)=0$, then $\lim\limits_{t\rightarrow+\infty}{I}(y,t)=0$
uniformly for $y\in\overline{\Omega(0)}$.

The above result implies that for any $\varepsilon>0$, there exists $T_\varepsilon>0$ such that $0\leq I(y,t)\leq
\varepsilon$ for $y\in \overline{\Omega(0)}$ and $t>T_\varepsilon$, then
{\small $$
-M^{*}\varepsilon+a(\rho(t)y)S-b(\rho(t)y)S^{2}-\frac{n\dot \rho(t)}{\rho(t)}S\leq S_{t}-\frac{d_S}{\rho^2(t)}
\Delta S\leq M^{*}\varepsilon+a(\rho(t)y)S-b(\rho(t)y)S^{2}-\frac{n\dot \rho(t)}{\rho(t)}S,
$$}
where $M^{*}=\max\limits_{(y,t)\in\overline{\Omega(0)}\times[0,T]}\{\beta(\rho(t)y)+\gamma(\rho(t)y)\}$.
Suppose that $\overline{S}_\varepsilon$ and $\underline{S}_\varepsilon$ satisfy the following problems
\begin{eqnarray*}
\left\{
\begin{array}{lllll}
S_{t}-\frac{d_S}{\rho^2(t)} \Delta S=M^{*}\varepsilon+a(\rho(t)y)S-b(\rho(t)y)S^{2}-\frac{n\dot \rho(t)}{\rho(t)}S,
\; &\ y\in \Omega(0), t>0, \\[2mm]
\frac{\partial S(y,t)}{\partial \nu}=0,\ & \ y\in \partial \Omega(0), t>0,\\[2mm]
S(y,0)=S_0(y)>0,\, \; &\ y\in  \overline{\Omega(0)},
\end{array} \right.
\end{eqnarray*}
and
\begin{eqnarray*}
\left\{
\begin{array}{lllll}
S_{t}-\frac{d_S}{\rho^2(t)} \Delta S=-M^{*}\varepsilon+a(\rho(t)y)S-b(\rho(t)y)S^{2}-\frac{n\dot \rho(t)}{\rho(t)}S,
\; &\ y\in \Omega(0), t>0, \\[2mm]
\frac{\partial S(y,t)}{\partial \nu}=0,\ & \ y\in \partial \Omega(0), t>0,\\[2mm]
S(y,0)=S_0(y)>0,\, \; &\ y\in \overline{\Omega(0)}.
\end{array} \right.
\end{eqnarray*}
Using comparison principle, we can derive that $\overline{S}_\varepsilon$ is the upper solution of problem
\eqref{a07}-\eqref{a08} and $\underline{S}_\varepsilon$ is the lower solution of problem \eqref{a07}-\eqref{a08}.
Hence, we can know that the solution $(S(y,t), I(y, t))$ of problem \eqref{a07}-\eqref{a08}
satisfies $\underline{S}_\varepsilon\leq S(y,t)\leq \overline{S}_\varepsilon$ in $\Omega(0)\times[0,+\infty)$
by the method of upper and lower solution. We obtain a sequence $\overline{S}_\varepsilon^{(m)}$ and
$\underline{S}_\varepsilon^{(m)}$ from the following the parabolic problem as initial iterations
$\overline{S}_\varepsilon^{(0)}=\overline{S}_\varepsilon$ and $\underline{S}_\varepsilon^{(0)}=\underline{S}_\varepsilon$
\begin{eqnarray}
\left\{
\begin{array}{lll}
(\overline{S}_\varepsilon^{(m)})_{t}-\frac{d_S}{\rho^2(t)} \Delta \overline{S}_\varepsilon^{(m)}+K\overline{S}_\varepsilon^{(m)}=g_1(\overline{S}_\varepsilon^{(m-1)}),
\; &\ y\in \Omega(0), t>0,\\[2mm]
(\underline{S}_\varepsilon^{(m)})_{t}-\frac{d_S}{\rho^2(t)} \Delta \underline{S}_\varepsilon^{(m)}+K\underline{S}_\varepsilon^{(m)}=g_2(\underline{S}_\varepsilon^{(m-1)}),
\; &\ y\in \Omega(0), t>0,\\[2mm]
\frac{\partial \overline{S}_\varepsilon^{(m)}(y,t)}{\partial \nu}=
\frac{\partial \underline{S}_\varepsilon^{(m)}(y,t)}{\partial \nu}=0,\; &\ y\in \partial\Omega(0), t>0,\\[2mm]
\overline{S}_\varepsilon^{(m)}(y,0)=\overline{S}_\varepsilon^{(m-1)}(y,T),\underline{S}_\varepsilon^{(m)}(y,0)
=\underline{S}_\varepsilon^{(m-1)}(y,T),\; &\ y\in \overline{\Omega(0)},
\end{array} \right.
\label{c09}
\end{eqnarray}
where $m=1,2,\cdots$ and
$$
g_1(S)=M^*\varepsilon+a(\rho(t)y)-b(\rho(t)y)S^{2}-\frac{n\dot \rho(t)}{\rho(t)}S+KS,
$$
$$
g_2(S)=-M^*\varepsilon+a(\rho(t)y)-b(\rho(t)y)S^{2}-\frac{n\dot \rho(t)}{\rho(t)}S+KS,
$$
$$
K=\sup\limits_{t\in[0,T)}\{\frac{n\dot \rho(t)}{\rho(t)}\}+2\max_{(y,t)\in(\overline{\Omega(0)}\times[0,+\infty))}
b(\rho(t)y)\overline{S}_{\varepsilon}
-\min_{(y,t)\in(\overline{\Omega(0)}\times[0,+\infty))}a(\rho(t)y).
$$

The sequences $\overline{S}_\varepsilon^{(m)}$ and $\underline{S}_\varepsilon^{(m)}$ satisfy the monotone property
$$
\underline{S}_\varepsilon\leq\underline{S}_\varepsilon^{(m-1)}\leq\underline{S}_\varepsilon^{(m)}
\leq\overline{S}_\varepsilon^{(m)}\leq\overline{S}_\varepsilon^{(m-1)}\leq\overline{S}_\varepsilon
$$
by Lemma 2.1 in \cite{pao}.
Due to the above inequality, it is easy to see that the pointwise limits
$$
\lim\limits_{m\rightarrow\infty}\overline{S}_\varepsilon^{(m)}=\overline{S}_\varepsilon^{*},\
\lim\limits_{m\rightarrow\infty}\underline{S}_\varepsilon^{(m)}=\underline{S}_\varepsilon^{*},
$$
which implies that
$$
\underline{S}_\varepsilon\leq\underline{S}_\varepsilon^{(m-1)}\leq\underline{S}_\varepsilon^{(m)}
\leq\underline{S}_\varepsilon^{*}\leq\overline{S}_\varepsilon^{*}\leq\overline{S}_\varepsilon^{(m)}
\leq\overline{S}_\varepsilon^{(m-1)}\leq\overline{S}_\varepsilon.
$$
Noticing that
\begin{equation*}
 \underline{S}_\varepsilon(y,t)\leq S(y,t) \leq \overline{S}_\varepsilon(y,t) \ \   \rm{in} \ \overline
 {\Omega(0)}\times [0,+\infty)
\end{equation*}
and setting $S_m(y,t)=S(y,t+mT)$, then
 $$
 \underline{S}_\varepsilon(y,t+T)\leq S(y,t+T)=S_1(y,t) \leq \overline{S}_\varepsilon(y,t+T) \ \   \rm{in} \
 \overline{\Omega(0)}\times [0,+\infty).
 $$
Replacing $S_{1}(y,0)$ with $S_{0}(y)$ in system \eqref{a07}, combining \eqref{c09} with $m=1$, then
$$
\overline{S}_\varepsilon^{(1)}(y, 0)=\overline{S}_\varepsilon^{(0)}(y, T)=\overline{S}_\varepsilon (y,T)
$$
and
$$
\underline{S}_\varepsilon^{(1)}(y,0)=\underline{S}_\varepsilon^{(0)}(y, T)=\underline{S}_\varepsilon(y,T).
$$
Combining
$$
\underline{S}_\varepsilon^{(1)}(y,0)=\underline{S}_\varepsilon(y,T)\leq S(y,T)=S_1(y,0)\leq\overline{S}_\varepsilon^{(1)}(y,0)=\overline{S}_\varepsilon(y,T)\   \rm{in} \ \Omega(0),
$$
with comparison principle obtain that
$$
\underline{S}_\varepsilon^{(1)}(y,t)\leq S_1(y,t)\leq\overline{S}_\varepsilon^{(1)}(y,t)\   \rm{in} \ \Omega(0)
\times [0,+\infty).
$$
Repeating the above iteration gives that
$$
\underline S_\varepsilon^{m}(y,t) \leq S_{m}(y,t) \leq \overline S_\varepsilon^{m}(y,t)\  \rm{in} \ \overline
{\Omega(0)}\times [0,+\infty).
$$
Consequently, we have
{\small $$
\underline{S}_\varepsilon^*(y,t)=\lim_{m\to\infty}\underline S_\varepsilon^{(m)}(y,t)\leq \liminf_{m\to\infty}
S_{m}(y,t)\leq \limsup_{m\to\infty} S_{m}(y,t) \leq \lim_{m\to\infty}\overline S_\varepsilon^{(m)}(y,t)=
\overline{S}_\varepsilon^*(y,t)
$$}
for $(y,t)\in \overline{\Omega(0)}\times [0,+\infty)$. Moreover, $\underline{S}_\varepsilon^*(y,t)$ satisfies
\begin{eqnarray*}
\left\{
\begin{array}{lll}
S_{t}-\frac{d_S}{\rho^2(t)} \Delta S=-M^*\varepsilon+a(\rho(t)y)S-b(\rho(t)y)S^{2}-\frac{n\dot \rho(t)}{\rho(t)}S,
\; &\ y\in \Omega(0), t>0,\\[2mm]
\frac{\partial S(y,t)}{\partial \nu}=0,\; &\ y\in \partial\Omega(0), t>0,\\[2mm]
S(y,0)=S(y,T),\; &\ y\in \overline{\Omega(0)},
\end{array} \right.
\end{eqnarray*}
and $\overline{S}_\varepsilon^*(y,t)$ satisfies
\begin{eqnarray*}
\left\{
\begin{array}{lll}
S_{t}-\frac{d_S}{\rho^2(t)} \Delta S=M^*\varepsilon+a(\rho(t)y)S-b(\rho(t)y)S^{2}-\frac{n\dot \rho(t)}{\rho(t)}S,
\; &\ y\in \Omega(0), t>0,\\[2mm]
\frac{\partial S(y,t)}{\partial \nu}=0,\; &\ y\in \partial\Omega(0), t>0,\\[2mm]
S(y,0)=S(y,T),\; &\ y\in \overline{\Omega(0)}.
\end{array} \right.
\end{eqnarray*}

According to the uniqueness of the solution to problem \eqref{b01}, we have
$$
\lim\limits_{\varepsilon \to 0^{+} } \overline{S}_\varepsilon^*(y,t)=\lim\limits_{\varepsilon \to 0^{+}}
\underline{S}_\varepsilon^*(y,t)=S^{*}(y,t),
$$
therefore
$$
\lim\limits_{m \to \infty } S(y,t+mT)=S^{*}(y,t) \ \textrm{for} \ \overline{\Omega(0)}\times [0,+\infty).
$$

$(ii)$ Due to $\lim\limits_{I\rightarrow0}f(S,I)=1$,
take $\delta=\frac{1}{4}(1-\frac{1}{{\mathcal{R}}_0})>0$, there exists $\varepsilon_0>0$ such that
$$
1-\delta\leq f(S,I)\leq 1
$$
if $ 0\leq I(y,t)\leq\varepsilon_0$.

We prove (ii) of Theorem \ref{main} by contradiction. Assume that there exists a nonnegative nontrivial
solution $(S, I)$ of problem \eqref{a07}, \eqref{a09} such that
\begin{equation}
\limsup\limits_{t\to\infty}\|(S(y,t),I(y,t))-(S^*(y,t),0)\|<\varepsilon_0/2.
\label{c10}
\end{equation}
For the above given $\varepsilon_0$, there exists $T_{\varepsilon_0}$ such that
$$
0\leq I(y,t)\leq\varepsilon_0\; \ \textrm{for} \;(y,t)\in\overline{\Omega(0)}\times [T_{\varepsilon_0}, \infty).
$$
Therefore, we have
\begin{equation}
\begin{array}{llllll}
I_{t}-\frac{d_I}{\rho^2(t)} \Delta I
&=&\beta(\rho(t)y)f(S,{I})I-\gamma(\rho(t)y)I-\frac{n\dot \rho(t)}{\rho(t)}I \\[2mm]
&\geq&\beta(\rho(t)y)(1-\delta)I-\gamma(\rho(t)y)I-\frac{n\dot \rho(t)}{\rho(t)}I
\end{array}
\label{c11}
\end{equation}
for $y\in\Omega(0),\; t\geq T_{\varepsilon_0}$. We now choose a small enough number $\alpha>0$ and
$\Phi(y,t)>0$ for $(y,t)\in\overline{\Omega(0)}\times[0,T]$ satisfies \eqref{c01} with ${\mathcal{R}}_0>1$
such that
\begin{equation}
I(y, T_{\varepsilon_0})\geq \alpha\Phi(y, T_{\varepsilon_0}).
\label{c12}
\end{equation}
Let $0<\lambda_0\leq\frac{3}{4}(1-\frac{1}{{\mathcal{R}}_0})\beta(\rho(t)y)$ and there exists a sufficiently
small constant $\eta>0$ such that $\eta e^{\lambda_0 T_{\varepsilon_0}}\leq\alpha$, then it is easy to show
that $\underline{I}(y,t)=\eta e^{\lambda_0t}\Phi(y,t)$ satisfies
\begin{eqnarray*}
\left\{
\begin{array}{lll}
\underline I_{t}-\frac{d_I}{\rho^2(t)} \Delta \underline I\leq \beta(\rho(t)y)(1-\delta)\underline I-
\gamma(\rho(t)y)\underline I-\frac{n\dot \rho(t)}{\rho(t)}\underline I,&\ y\in \Omega(0), t\geq T_{\varepsilon_0}, \\[2mm]
\frac{\partial \underline I(y,t)}{\partial \nu}=0,\ & \ y\in \partial \Omega(0), t\geq T_{\varepsilon_0},\\[2mm]
\underline I(y, T_{\varepsilon_0}) \leq I(y, T_{\varepsilon_0}), \; &\ y\in  \overline{\Omega(0)}.
\end{array} \right.
\end{eqnarray*}
In light of \eqref{c11}, \eqref{c12} and comparison principle, we get
$$
I(y,t)\geq \underline I(y, t)=\eta e^{\lambda_0t}\Phi(y,t) \;\rm{for}\;   \textit{y}\in\Omega(0),\
\textit{t}\geq T_{\varepsilon_0},
$$
therefore, $I(y,t) \rightarrow\infty$ as $t\rightarrow\infty$, which contradicts \eqref{c10}.
This proves statement $(ii)$.
\epf

\section{\bf The effects of $d_I$, $L$ and $\rho(t)$ on ${\mathcal{R}}_0$}
This section is devoted to researching the effects of diffusion coefficient $d_I$ and interval length $L$ on
${\mathcal{R}}_0$ in one-dimensional space, and evolving rate $\rho(t)$ on ${\mathcal{R}}_0$ in some special cases.
We first assume that $\gamma_{y}(\rho(t)y)<0$ and $\beta_{y}(\rho(t)y)\geq0$ for $(y,t)\in\overline{\Omega(0)}\times[0,T]$.
Problem \eqref{c01} in one-dimensional space becomes
\begin{eqnarray}
\left\{
\begin{array}{ll}
\Phi_{t}-\frac{d_I}{\rho^2(t)} \Phi_{yy}+[\gamma(\rho(t)y)+\frac{\dot \rho(t)}{\rho(t)}]\Phi=\frac{\beta(\rho(t)y)}
{{\mathcal{R}}_0}\Phi,\; &\ 0<y<L, 0<t\leq T, \\[2mm]
\frac{\partial \Phi(y,t)}{\partial \nu}=0,\ & \ y=0,L, 0<t\leq T,\\[2mm]
\Phi(y,0)=\Phi(y,T),\; &\ 0\leq y\leq L.
\end{array} \right.
\label{d01}
\end{eqnarray}

According to Lemma 1 of \cite{pz1}, we can obtain following result.
\begin{prop}
\label{case-2.2} $\Phi_{y}(y,t)>0$ in $(0,L)$ for all $t$ and for any given $d_I, L>0$.
\end{prop}

Besides, we consider the adjoint problem of \eqref{d01}
\begin{eqnarray}
\left\{
\begin{array}{ll}
-\phi_{t}-\frac{d_I}{\rho^2(t)}\phi_{yy}+[\gamma(\rho(t)y)+\frac{\dot \rho(t)}{\rho(t)}]\phi=\frac{\beta(\rho(t)y)}
{{\mathcal{R}}_0}\phi,\; &\ 0<y<L, 0<t\leq T, \\[2mm]
\frac{\partial \phi(y,t)}{\partial \nu}=0,\ & \ y=0,L, 0<t\leq T,\\[2mm]
\phi(y,0)=\phi(y,T),\; &\ 0\leq y\leq L.
\end{array} \right.
\label{d02}
\end{eqnarray}
Therefore, problems \eqref{d01} and \eqref{d02} have the same principal eigenvalue ${\mathcal{R}}_0$.

We now investigate the monotonicity of ${\mathcal{R}}_0$ in problem \eqref{d01} with respect to $d_I$ and $L$.
\begin{thm}
\label{main1}
The following statements hold:

$(i)$ ${\mathcal{R}}_0$ is strictly monotone decreasing in $d_I$ for any given $L>0$;

$(ii)$ ${\mathcal{R}}_0$ is strictly monotone increasing in $L$ for any given $d_I>0$.

\end{thm}
\bpf  We first prove $(i)$. Obviously, ${\mathcal{R}}_0$ and $\Phi$ are $C^1$-functions of $d_I$. For the
sake of simplicity, we use $\Phi'$ denote $\frac{\partial \Phi(y,t)}{\partial d_I}$ and ${\mathcal{R}}_0'$
denote $\frac{\partial {\mathcal{R}}_0}{\partial d_I}$. Differentiating \eqref{d01} with respect to $d_I$ yields
\begin{small}
\begin{eqnarray}
\left\{
\begin{array}{ll}
\Phi_{t}'-\frac{d_I}{\rho^2(t)} \Phi'_{yy}-\frac{1}{\rho^2(t)}\Phi_{yy}=\frac{\beta(\rho(t)y)}{{\mathcal{R}}_0}\Phi'+(\frac{\beta(\rho(t)y)}{{\mathcal{R}}_0})'\Phi, & \\
\qquad \qquad \qquad \qquad \quad -[\gamma(\rho(t)y)+\frac{\dot \rho(t)}{\rho(t)}]\Phi'  \; &\ 0<y<L, 0<t\leq T, \\[2mm]
\frac{\partial \Phi'(y,t)}{\partial \nu}=0,\ & \ y=0,L, 0<t\leq T,\\[2mm]
\Phi'(y,0)=\Phi'(y,T),\; &\ 0\leq y\leq L.
\end{array} \right.
\label{cth101}
\end{eqnarray}
\end{small}
Multiplying the equation of \eqref{cth101} by $\phi$, which satisfies \eqref{d02} and then integrating
the resulting equation give
$$\begin{array}{lll}
&&-\int_0^{T}\int_0^{L}\Phi'\phi_{t} dydt-\int_0^{T}\int_0^{L}\frac{\Phi_{yy}}{\rho^2(t)}\phi dydt\\[2mm]
&-&\int_0^{T}\int_0^{L}\frac{d_I}{\rho^2(t)}\Phi'_{yy}\phi dydt+\int_0^{T}\int_0^{L}[\gamma(\rho(t)y)+
\frac{\dot \rho(t)}{\rho(t)}]\Phi'\phi dydt\\[2mm]
&=&\int_0^{T}\int_0^{L}\frac{1}{{\mathcal{R}}_0}\beta(\rho(t)y)\Phi'\phi dydt+
\int_0^{T}\int_0^{L}(\frac{1}{{\mathcal{R}}_0})'\beta(\rho(t)y)\Phi\phi dydt.
\end{array}$$
Substituting $-\phi_{t}-\frac{d_I}{\rho^2(t)}\phi_{yy}+[\gamma(\rho(t)y)+
\frac{\dot \rho(t)}{\rho(t)}]\phi=\frac{\beta(\rho(t)y)}{{\mathcal{R}}_0}\phi$ into the above equation yields
$$
(\frac{1}{{\mathcal{R}}_0})'=\frac{-\int_0^{T}\int_0^{L}\frac{\Phi_{yy}}{\rho^2(t)}\phi dydt}
{\int_0^{T}\int_0^{L}\beta(\rho(t)y)\Phi\phi dydt}=\frac{\int_0^{T}\int_0^{L}\frac{\Phi_{y}\phi_{y}}{\rho^2(t)}dydt}{\int_0^{T}\int_0^{L}\beta(\rho(t)y)\Phi\phi dydt}.
$$
It follows proposition \ref{case-2.2} that $\Phi_{y}(y,t)>0$. To derive our result  ${\mathcal{R}}_0'<0$ in $(i)$, it suffices to prove that $\phi_{y}>0$.

In fact, by \eqref{d02}, let $\vartheta(y,t)=\phi(y,T-t)$, then $\vartheta(y,t)$
satisfies
\begin{equation*}
\left\{
\begin{array}{ll}
\vartheta_{t}-\frac{d_I}{\rho^2(T-t)} \vartheta_{yy}+[\gamma(\rho(T-t)y)+\frac{\dot \rho(T-t)}{\rho(T-t)}]\vartheta=\frac{\beta(\rho(T-t)y)}{{\mathcal{R}}_0}\vartheta,\; &\ 0<y<L, 0<t\leq T, \\[2mm]
\frac{\partial \vartheta(y,t)}{\partial \nu}=0,\ & \ y=0,L, 0<t\leq T,\\[2mm]
\vartheta(y,0)=\vartheta(y,T),\; &\ 0\leq y\leq L.
\end{array} \right.
\end{equation*}
Due to $\gamma_{y}(\rho(t)y)<0$ and $\beta_{y}(\rho(t)y)\geq0$ for $(y,t)\in[0,L]\times[0,T]$ and Proposition
\ref{case-2.2}, we can get
$$
\vartheta_y=\phi_y(y,T-t)>0,\ \ \  \forall (y,t)\in (0,L)\times[0,T].
$$
The prove of $(i)$ is completed.

We now prove $(ii)$. Suppose that $0<L_1<L_2$, $\Phi$ is the principal eigenfunction of \eqref{d01} corresponding to
${\mathcal{R}}_0(L_2)$, and $\phi$ is the principal eigenfunction of \eqref{d02} corresponding to
${\mathcal{R}}_0(L_1)$. Hence, $(\Phi,L_2)$ solves
\begin{equation*}
\left\{
\begin{array}{ll}
\Phi_{t}-\frac{d_I}{\rho^2(t)} \Phi_{yy}+[\gamma(\rho(t)y)+\frac{\dot \rho(t)}{\rho(t)}]\Phi=\frac{\beta(\rho(t)y)}
{{\mathcal{R}}_0(L_2)}\Phi,\; &\ 0<y<L_2, 0<t\leq T, \\[2mm]
\frac{\partial \Phi(y,t)}{\partial \nu}=0,\ & \ y=0,L_2, 0<t\leq T,\\[2mm]
\Phi(y,0)=\Phi(y,T),\; &\ 0\leq y\leq L_2,
\end{array} \right.
\end{equation*}
and $(\phi,L_1)$ satisfies
\begin{equation*}
\left\{
\begin{array}{ll}
-\phi_{t}-\frac{d_I}{\rho^2(t)}\phi_{yy}+[\gamma(\rho(t)y)+\frac{\dot \rho(t)}{\rho(t)}]\phi=\frac{\beta(\rho(t)y)}
{{\mathcal{R}}_0(L_1)}\phi,\; &\ 0<y<L_1, 0<t\leq T, \\[2mm]
\frac{\partial \phi(y,t)}{\partial \nu}=0,\ & \ y=0,L_1, 0<t\leq T,\\[2mm]
\phi(y,0)=\phi(y,T),\; &\ 0\leq y\leq L_1.
\end{array} \right.
\end{equation*}
We multiply the equation of $\Phi$ by $\phi$, the equation of $\phi$ by $\Phi$, and subtract the resulting equations
to obtain
$$
(\frac{1}{{\mathcal{R}}_0(L_2)}-\frac{1}{{\mathcal{R}}_0(L_1)})\beta(\rho(t)y)\Phi\phi=\Phi_t\phi+\phi_t\Phi+
\frac{d_I}{\rho^{2}(t)}\phi_{yy}\Phi-\frac{d_I}{\rho^{2}(t)}\Phi_{yy}\phi
$$
in $(0,L_1)\times(0,T)$.
Integrating the above equality over $(0,L_1)\times(0,T)$ yields
$$\begin{array}{lll}
&&(\frac{1}{{\mathcal{R}}_0(L_2)}-\frac{1}{{\mathcal{R}}_0(L_1)})\int_0^{T}\int_0^{L_1}\beta(\rho(t)y)\Phi\phi dydt\\[2mm]
&=&\int_0^{T}\int_0^{L_1}\frac{d_{I}}{\rho^2(t)}(\phi_{yy}\Phi- \Phi_{yy}\phi)dydt+\int_0^{T}\int_0^{L_1}(\Phi_t\phi
+\phi_t\Phi)dydt\\[2mm]
&=&\int_0^{T}\frac{d_I}{\rho^2(t)}[\phi_y(L_1)\Phi(L_1)-\phi_y(0)\Phi(0)-\Phi_y(L_1)\phi(L_1)+\Phi_y(0)\phi(0) ]dt\\[2mm]
&=&-\int_0^{T}\frac{d_I}{\rho^2(t)}\Phi_y(L_1)\phi(L_1)dt.
\end{array}$$
Based on $\Phi_y(y,t)>0$ for all $(y,t)\in(0,L_1]\times[0,T]$, we can know
$$
\frac{1}{{\mathcal{R}}_0(L_2)}-\frac{1}{{\mathcal{R}}_0(L_1)}=
-\frac{\int_0^{T}\frac{d_I}{\rho^2(t)}\Phi_{y}(L_1)\phi(L_1)dt}{\int_0^{T}\int_0^{L_1}\beta(\rho(t)y)\Phi\phi dydt}<0.
$$
Therefore ${\mathcal{R}}_0(L)$ is strictly monotone increasing in $L$.

\epf
\medskip

In the following, we discuss the limiting behaviors of ${\mathcal{R}}_0$ as $d_I\rightarrow0$ or $d_I\rightarrow\infty$.

\begin{thm}
\label{main2}
The following statements hold:

$(i)$  $\lim\limits_{d_I\to 0}{{\mathcal{R}}_0}=\frac{\int^T_0\beta(\rho(t)L)dt}{\int^T_0\gamma(\rho(t)L)dt}$;

$(ii)$ $\lim\limits_{d_I\to\infty}{{\mathcal{R}}_0}=\frac{\int_0^{L}\int_0^{T}\beta(\rho(t)y)dtdy}
{\int_0^{L}\int_0^{T}\gamma(\rho(t)y)dtdy}$.

\end{thm}

\bpf
We first prove $(i)$. Consider the following auxiliary system
\begin{eqnarray}
\left\{
\begin{array}{ll}
\Phi_{t}-\frac{d_I}{\rho^2(t)} \Phi_{yy}=[\frac{\beta(\rho(t)y)}{{\mathcal{R}}_0}-\gamma(\rho(t)y)-
\frac{\dot \rho(t)}{\rho(t)}]\Phi+\Lambda\Phi,\; &\ 0<y<L, 0<t\leq T, \\[2mm]
\frac{\partial \Phi(y,t)}{\partial \nu}=0,\ & \ y=0,L, 0<t\leq T,\\[2mm]
\Phi(y,0)=\Phi(y,T),\; &\ 0\leq y\leq L.
\end{array} \right.
\label{cth2-1}
\end{eqnarray}
Assume that $\Lambda_1({\mathcal{R}}_0,d_I)$ is the principal eigenvalue of \eqref{cth2-1} and
$\Lambda^D({\mathcal{R}}_0,d_I)$ is the principal eigenvalue of \eqref{cth2-1} with Dirichlet
boundary condition. It is obvious that $\Lambda_1({\mathcal{R}}_0,d_I)\leq\Lambda^D({\mathcal{R}}_0,d_I)$
by the eigenvalue property. We say that $\mu (G(t))$ is an eigenvalue if there is a nontrivial
solution $\zeta$ to the following problem
\begin{eqnarray}
\left\{
\begin{array}{ll}
\zeta_{t}=G(t)\zeta+\mu\zeta,\; &\  0<t\leq T, \\[2mm]
\zeta(0)=\zeta(T).
\end{array} \right.
\label{cth2-2}
\end{eqnarray}
We notice that equations in \eqref{cth2-1} become equations in \eqref{d01} if $\Lambda({\mathcal{R}}_0,d_I)=0$.
Define $\mu_1(y):=\mu(M(y,\cdot))=-\frac{1}{T}\int_0^TM(y,t)dt$ is an eigenvalue of \eqref{cth2-2} with $G(t)$
replaced by $M(y,t)$, where $M(y,t)=\frac{\beta(\rho(t)y)}{{\mathcal{R}}_0}-\gamma(\rho(t)y)-\frac{\dot \rho(t)}{\rho(t)}$.
Our proof divides two steps.

\textbf{Step 1.} $\limsup\limits_{d_I\to0}\Lambda_1({\mathcal{R}}_0,d_I)\leq-\frac{1}{T}\int^T_0\max\limits_{[0,L]}M(y,t)dt$.
We now assume that $\min\limits_{[0,L]}\mu_1(y)=-\frac{1}{T}\int^T_0\max\limits_{[0,L]}M(y,t)dt<0$. By Lemma 3.6 in \cite{BH},
there exists $y_0\in(0,L)$ and $\delta>0$ small such that
\begin{equation}
\mu_1(y)<0  \ \ \ \ \forall y\in \overline{U_\delta(y_0)}\subset(0,L).
\label{cth2-3}
\end{equation}
Consider an auxiliary problem
\begin{eqnarray}
\left\{
\begin{array}{ll}
\Phi_{t}-\frac{d_I}{\rho^2(t)} \Phi_{yy}=M(y,t)\Phi+\widehat{\Lambda}\Phi,\; &\ y\in U_\delta(y_0), 0<t\leq T, \\[2mm]
\Phi(y,t)=0,\ & \ y\in \partial U_\delta(y_0), 0<t\leq T,\\[2mm]
\Phi(y,0)=\Phi(y,T),\; &\ y\in U_\delta(y_0).
\end{array} \right.
\label{cth2-4}
\end{eqnarray}
Denote $\widehat{\Lambda}_1({\mathcal{R}}_0,d_I)$ is the principal eigenvalue of the above problem. It
is obvious that $\Lambda^D({\mathcal{R}}_0,d_I)\leq\widehat{\Lambda}_1({\mathcal{R}}_0,d_I)$ by Corollary
4.2 in \cite{BH}. In light of \eqref{cth2-3}, it is easy to see that
$$
\mu_{1}(y_0)<0,
$$
where $\mu_{1}(y_0)$ is the principal eigenvalue for $M(y_0,\cdot)$. Let $\xi(y_0)$ is the positive principal
eigenfunction corresponding to $\mu_{1}(y_0)$.
Now choosing $\delta$ small enough and assuming
$$
\|M(y_0,t)-M(y,t)\|_{\overline{U}_{\delta(y_{0})}\times [0,T]}\leq-\frac{\mu_{1}(y_0)}{2}.
$$
There exists some $\widetilde{\psi}\in C^2(\overline{U}_{\delta(y_{0})})$ which is positive in
$U_{\delta(y_{0})}$ with $\widetilde{\psi}=0$ on $\partial U_{\delta(y_{0})}$,
and a positive constant $C$ such that
$$
\widetilde{\psi}_t-\frac{d_I}{\rho^2(t)}\widetilde{\psi}_{yy}=-\frac{d_I}{\rho^2(t)}\widetilde{\psi}_{yy}\leq
d_IC\widetilde{\psi}
$$
by Corollary 2.7 in \cite{BH}.
Directing calculation yields
$$
\begin{array}{llllll}
&&(\widetilde{\psi}\xi(y_0))_{t}-\frac{d_I}{\rho^2(t)} (\widetilde{\psi}\xi(y_0))_{yy}\\[2mm]
&=&\xi(y_0)_{t}\widetilde{\psi}+\widetilde{\psi}_{t}\xi(y_0)-\xi(y_0)\frac{d_I}{\rho^2(t)} \widetilde{\psi}_{yy}\\[2mm]
&\leq&M(y,t)\widetilde{\psi}\xi(y_0)+(M(y_0,t)-M(y,t))\widetilde{\psi}\xi(y_0)+\mu_1(y_0)\widetilde{\psi}\xi(y_0)
+d_IC\widetilde{\psi}\xi(y_0)\\[2mm]
&\leq&M(y,t)\widetilde{\psi}\xi(y_0)+\frac{\mu_1(y_0)}{2}\widetilde{\psi}\xi(y_0)+d_IC\widetilde{\psi}\xi(y_0)\\[2mm]
&\leq& M(y,t)\widetilde{\psi}\xi(y_0),
\end{array}
$$
where $\mu_1(y_0)<0$ and $d_I$ is sufficiently small such that $C\leq\frac{-\mu_1(y_0)}{d_I}$.
Hence, we can obtain
\begin{eqnarray}
\left\{
\begin{array}{ll}
(\widetilde{\psi}\xi(y_0))_{t}-\frac{d_I}{\rho^2(t)} (\widetilde{\psi}\xi(y_0))_{yy}-M\widetilde{\psi}\xi(y_0)\leq0,
\; &\ y\in U_\delta(y_0), 0<t\leq T, \\[2mm]
\widetilde{\psi}\xi(y_0)=0,\ & \ y\in \partial U_\delta(y_0), 0<t\leq T,\\[2mm]
(\widetilde{\psi}\xi(y_0))(y,0)=(\widetilde{\psi}\xi(y_0))(y,T),\; &\ y\in U_\delta(y_0).
\end{array} \right.
\label{cth2-5}
\end{eqnarray}
Combining \eqref{cth2-4} with \eqref{cth2-5}, and then applying Proposition 2.4 in \cite{BH} give
$$
\widehat{\Lambda}_1({\mathcal{R}}_0,d_I)\leq0.
$$
So
$$
\Lambda_1({\mathcal{R}}_0,d_I)\leq\Lambda^D({\mathcal{R}}_0,d_I)\leq\widehat{\Lambda}_1({\mathcal{R}}_0,d_I)\leq0.
$$
Let $\varepsilon>0$ and substitute $\mu^{\varepsilon}_{1}(y)=\mu_1(y)-\min\limits_{[0,L]}\mu_1(y)-\varepsilon$ for
$\mu_1(y)$ and substitute $M^\varepsilon(y,t)=M(y,t)+\min\limits_{[0,L]}\mu_1(y)+\varepsilon$ for $M(y,t)$. It is
easily seen that $\min\limits_{[0,L]}\mu^{\varepsilon}_1(y)=-\varepsilon<0$ and $\Lambda^{\varepsilon,D}_1({\mathcal{R}}_0,d_I)=\Lambda^D({\mathcal{R}}_0,d_I)-\min\limits_{[0,L]}\mu_1(y)-\varepsilon<0$,
where $\Lambda^{\varepsilon,D}_1({\mathcal{R}}_0,d_I)$ is the principal eigenvalue of \eqref{cth2-1} with Dirichlet
boundary condition and $M$ replaced by $M^\varepsilon$. If $d_I$ is small enough, then $\Lambda^D({\mathcal{R}}_0,d_I)\leq\min\limits_{[0,L]}\mu_1(y)+\varepsilon$
for any $\varepsilon>0$. Therefore,
$$
\limsup\limits_{d_I\to0}\Lambda_1({\mathcal{R}}_0,d_I)\leq\limsup\limits_{d_I\to0}\Lambda^D({\mathcal{R}}_0,d_I)
\leq\min\limits_{[0,L]}\mu_1(y)
$$
as $d_I\rightarrow0$ and $\varepsilon\rightarrow0$.

\textbf{Step 2.} $\liminf\limits_{d_I\to0}\Lambda_1({\mathcal{R}}_0,d_I)\geq-\frac{1}{T}\int^T_0\max\limits_{[0,L]}M(y,t)dt$.

Note that
$$
\Phi_{t}-\frac{d_I}{\rho^2(t)}\Phi_{yy}=M(y,t)\Phi+\Lambda\Phi\leq\Phi\max\limits_{[0,L]}M(y,t)+\Lambda\Phi.
$$
We consider following auxiliary eigenvalue problem
\begin{eqnarray}
\left\{
\begin{array}{ll}
\overline{\Phi}_{t}-\frac{d_I}{\rho^2(t)} \overline{\Phi}_{yy}=\max\limits_{[0,L]}M(y,t)\overline{\Phi}+
\widetilde{\Lambda}\overline{\Phi},\; &\ 0<y<L, 0<t\leq T, \\[2mm]
\frac{\partial \overline{\Phi}(y,t)}{\partial \nu}=0,\ & \ y=0,L, 0<t\leq T,\\[2mm]
\overline{\Phi}(y,0)=\overline{\Phi}(y,T),\; &\ 0\leq y\leq L,
\end{array} \right.
\label{cth2-6}
\end{eqnarray}
Let $\widetilde{\Lambda}_1({\mathcal{R}}_0,d_I)$ is principal eigenvalue of above problem.
Straightforward calculation gives
$$
\widetilde{\Lambda}_1({\mathcal{R}}_0,d_I)=-\frac{1}{T}\int^T_0\max\limits_{[0,L]}M(y,t)dt.
$$
Using Corollary 4.2 in \cite{BH} yields
$$
\Lambda_1({\mathcal{R}}_0,d_I)\geq\widetilde{\Lambda}_1({\mathcal{R}}_0,d_I)=
-\frac{1}{T}\int^T_0\max\limits_{[0,L]}M(y,t)dt.
$$
Due to above inequality, it is easy to see that
$$
\liminf\limits_{d_I\to0}\Lambda_1({\mathcal{R}}_0,d_I)\geq\liminf\limits_{d_I\to0}
\widetilde{\Lambda}_1({\mathcal{R}}_0,d_I)
=-\frac{1}{T}\int^T_0\max\limits_{[0,L]}M(y,t)dt.
$$

Due to Step 1 and Step 2, we know that
$$
\lim\limits_{d_I\to0}\Lambda_1({\mathcal{R}}_0,d_I)=-\frac{1}{T}\int^T_0\max\limits_{[0,L]}M(y,t)dt.
$$
By the monotonicity of $\gamma$ and $\beta$, we derive that
$$
\begin{array}{llllll}
&&-\frac{1}{T}\int^T_0\max\limits_{[0,L]}M(y,t)dt\\[2mm]
&=&-\frac{1}{T}\int^T_0\max\limits_{[0,L]}[\frac{\beta(\rho(t)y)}{{\mathcal{R}}_0}
-\gamma(\rho(t)y)-\frac{\dot \rho(t)}{\rho(t)}]dt\\[2mm]
&=&\frac{1}{T}\int^T_0\gamma(\rho(t)L)dt -\frac{\frac{1}{T}\int^T_0\beta(\rho(t)L)dt}{{\mathcal{R}}_0}.
\end{array}
$$
For any $\varepsilon>0$, there exists a $\delta>0$ such that
$$
|\Lambda_1({\mathcal{R}}_0,d_I)-\frac{1}{T}\int^T_0\gamma(\rho(t)L)dt +\frac{\frac{1}{T}\int^T_0\beta(\rho(t)L)dt}{{\mathcal{R}}_0}|<\varepsilon
$$
for $0<d_I<\delta$.
It is easily seen that
$$
\frac{\frac{1}{T}\int^T_0\gamma(\rho(t)L)dt-\varepsilon-\Lambda_1({\mathcal{R}}_0,d_I)}{\frac{1}{T}\int^T_0\beta(\rho(t)L)dt}<{\mathcal{R}}_0<
\frac{\frac{1}{T}\int^T_0\gamma(\rho(t)L)dt+
\varepsilon-\Lambda_1({\mathcal{R}}_0,d_I)}{\frac{1}{T}\int^T_0\beta(\rho(t)L)dt}.
$$
Letting $\varepsilon\rightarrow0$ and $d_I\rightarrow0$ give
$$
\lim\limits_{d_I\to0}{\mathcal{R}}_0=\frac{\int^T_0\beta(\rho(t)L)dt}{\int^T_0\gamma(\rho(t)L)dt}
$$
when $\Lambda_1({\mathcal{R}}_0,d_I)=0$.

We then prove $(ii)$. Define
$$
\rho^{M}=\max\limits_{t\in[0,T]}\rho(t).
$$
At first, we normalize $\Phi$ such that
\begin{equation}
\int_0^T\int_0^L\Phi^{2}dydt=1.
\label{cth201}
\end{equation}
Multiplying the equation of \eqref{d01} by $\Phi$ and integrating over $(0,L)\times(0,T)$ give
$$
\int_0^T\int_0^L\frac{d_I}{\rho^{2}(t)}\Phi^{2}_{y}dydt=\int_0^T\int_0^L[\frac{\beta(\rho(t)y)}{{\mathcal{R}}_0}
-\gamma(\rho(t)y)-\frac{\dot{\rho(t)}}{\rho(t)}]\Phi^{2}dydt.
$$
Using the boundedness of ${\mathcal{R}}_0$ in \eqref{c08} and \eqref{cth201} give
$$
\frac{d_I}{(\rho^{M})^2}\int_0^{T}\int_0^{L}\Phi^{2}_{y}dydt\leq\int_0^{T}\int_0^{L}\frac{d_I}{\rho^{2}(t)}\Phi^{2}_{y}dydt\leq C,
$$
i.e.,
\begin{equation}
\int_0^{T}\int_0^{L}\Phi^{2}_{y}dydt\leq\frac{C(\rho^{M})^2}{d_I},
\label{cth202}
\end{equation}
where $C$ is positive constant and independent of $d_I$.

In addition, let
$$
\Phi_{*}(t):=\frac{1}{L}\int_0^{L}\Phi(y,t)dy \ \ \mbox{and} \ \ \Phi^{*}(y,t):=\Phi(y,t)-\Phi_{*}(t).
$$
Obviously,
$$
\int_0^{L}\Phi^{*}(y,t)dy=0 \ \ \mbox{for all} \ \ t\in [0,T].
$$
Using the well-known Poincar$\acute{e}$ inequality we conclude that
$$
\int_0^{L}(\Phi^{*}(y,t))^{2}dy \leq C\int_0^{L}(\Phi^{*}_{y})^{2}dy \ \ \mbox{for all} \ \ t\in [0,T].
$$
Making use of \eqref{cth202} and $\Phi^{*}_{y}=\Phi_{y}$, we obtain
\begin{equation*}
\int_0^{T}\int_0^{L}(\Phi^{*})^{2}dydt\leq\frac{C(\rho^{M})^2}{d_I},
\end{equation*}
which together with \eqref{cth202} yields by H$\ddot{o}$lder inequality
\begin{equation}
\int_0^{T}\int_0^{L}|\Phi^{*}|dydt\leq C\rho^{M}(\frac{1}{d_I})^{\frac{1}{2}},\ \ \
\int_0^{T}\int_0^{L}|\Phi^{*}_{y}|dydt\leq C\rho^{M}(\frac{1}{d_I})^{\frac{1}{2}}.
\label{cth203}
\end{equation}
On the other hand, integrating the equation of \eqref{d01} over $(0,L)$, we have
\begin{equation}
L\frac{d\Phi_{*}}{dt}=\Phi_{*}\int_0^{L}[\frac{\beta(\rho(t)y)}{{\mathcal{R}}_0}-\gamma(\rho(t)y)-\frac{\dot \rho(t)}{\rho(t)}]dy
+\int_0^{L}[\frac{\beta(\rho(t)y)}{{\mathcal{R}}_0}-\gamma(\rho(t)y)-\frac{\dot \rho(t)}{\rho(t)}]\Phi^*dy.
\label{cth204}
\end{equation}

Solving the ODE \eqref{cth204} gives us
\begin{equation}
\Phi_{*}(t)=\Phi_{*}(0)\exp(\frac{1}{L}\int_0^t\int_0^L(\frac{\beta(\rho(t)y)}{{\mathcal{R}}_0}-\gamma(\rho(t)y)-\frac{\dot \rho(t)}{\rho(t)})dydt)+\eta(t),
\label{cth205}
\end{equation}
where
{\small $$
\eta(t)=\frac{1}{L}\int_0^t\int_0^L\exp(\frac{1}{L}\int_\tau^t\int_0^L(\frac{\beta(\rho(t)y)}{{\mathcal{R}}_0}
-\gamma(\rho(t)y)-\frac{\dot \rho(t)}{\rho(t)})dydt)[\frac{\beta(\rho(t)y)}{{\mathcal{R}}_0}-\gamma(\rho(t)y)
-\frac{\dot \rho(t)}{\rho(t)}]\Phi^*dydt
$$}
for $\tau,t\in(0,T)$.

Due to \eqref{c08}, there exists a constant $M>0$ such that
$$
\exp (\frac{1}{L}\int_\tau^t\int_0^L(\frac{\beta(\rho(t)y)}{{\mathcal{R}}_0}-\gamma(\rho(t)y)-
\frac{\dot \rho(t)}{\rho(t)})dydt)\leq M.
$$
Using \eqref{cth203} and above inequality, it is easy to see that
$$
\eta(t)\leq M\int_0^T\int_0^L|[\frac{\beta(\rho(t)y)}{{\mathcal{R}}_0}-\gamma(\rho(t)y)-
\frac{\dot \rho(t)}{\rho(t)}]\Phi^*|dydt
\rightarrow0 \ \mbox{uniformly on}\, [0,T] \  \mbox{as} \ d_I\rightarrow\infty.
$$

Choosing $t=T$ in \eqref{cth205} and using $\Phi_{*}(0)=\Phi_{*}(T)$, we know
$$
\Phi_{*}(0)[1-\mbox{exp}(\frac{1}{L}\int_0^T\int_0^L(\frac{\beta(\rho(t)y)}{{\mathcal{R}}_0}-
\gamma(\rho(t)y)-\frac{\dot \rho(t)}{\rho(t)})dydt)]=\eta(T).
$$
Hence, when $d_I\rightarrow\infty$, then
$$
\mbox{either}\ \ \ \Phi_{*}(0)\rightarrow0 \ \ \mbox{or} \ \ \int_0^T\int_0^L[\frac{\beta(\rho(t)y)}{{\mathcal{R}}_0}
-\gamma(\rho(t)y)-\frac{\dot \rho(t)}{\rho(t)}]dydt\rightarrow 0.
$$
Suppose $\Phi_{*}(0)\rightarrow0$ as $d_I\rightarrow\infty$. In view of \eqref{cth205}, we easily see that
\begin{equation}
\Phi_{*}(t)\rightarrow0 \ \ \mbox{uniformly on}\ \  [0,T] \ \ \mbox{as} \ \ d_I\rightarrow\infty.
\label{cth206}
\end{equation}
Combining \eqref{cth203} with \eqref{cth206}, ensures
$$
\int_0^T\int_0^L\Phi dydt\rightarrow 0 \ \ \mbox{as} \ \ \ d_I\rightarrow \infty,
$$
which is a contradiction with \eqref{cth201}.

Hence,
$$
\int_0^T\int_0^L[\frac{\beta(\rho(t)y)}{{\mathcal{R}}_0}-\gamma(\rho(t)y)-\frac{\dot \rho(t)}{\rho(t)}]dydt\rightarrow0
\ \ \mbox{as}\ d_I\rightarrow\infty.
$$
i.e.,
$$
{\mathcal{R}}_0\rightarrow\frac{\int_0^T\int_0^L\beta(\rho(t)y)}{\int_0^T\int_0^L\gamma(\rho(t)y)}\ \
\mbox{as}\ d_I\rightarrow\infty.
$$

\epf
\medskip

From now on, we discuss the effect of $L$ on the basic reproduction number ${\mathcal{R}}_0$.

\begin{thm}
\label{main3}
The following statements hold:

$(i)$ $\lim\limits_{L\to0}{{\mathcal{R}}_0}=\frac{\int_0^{T}\beta(0)dt}{\int_0^{T}\gamma(0)dt}$.

$(ii)$ Assume $\lim\limits_{z\to \infty}\beta(z)=\beta_{\infty}$ and $\lim\limits_{z\to \infty}\gamma(z)
=\gamma_{\infty}$,
then $\lim\limits_{L\to \infty}{{\mathcal{R}}_0}=\frac{\beta_{\infty}}{\gamma_{\infty}}$.
\end{thm}
\bpf
We first prove $(i)$. Assuming that $({\mathcal{R}}_0,\Phi)$ satisfies the equation \eqref{d01}. Due to
the monotonicity of $\beta$ and $\gamma$, for any $\epsilon>0$, let $0<L<\epsilon$ to be sufficiently
small such that
$$
\beta(0)\leq\beta(\rho(t)y)\leq\beta(\rho(t)\epsilon) \ \mbox{on}\ [0,L]\times[0,T]
$$
and
$$
\gamma(\rho(t)\epsilon)\leq\gamma(\rho(t)y)\leq\gamma(0) \ \mbox{on}\ [0,L]\times[0,T].
$$
Thus $({\mathcal{R}}_0,\Phi)$ solves
\begin{equation*}
 \Phi_{t}-\frac{d_I}{\rho^2(t)} \Phi_{yy}+[\gamma(0)+\frac{\dot \rho(t)}{\rho(t)}]\Phi\geq\frac{\beta(0)}{{\mathcal{R}}_0}\Phi,
\ \ \   (y,t)\in[0,L]\times[0,T]
\end{equation*}
and
$$
\Phi_{t}-\frac{d_I}{\rho^2(t)} \Phi_{yy}+[\gamma(\rho(t)\epsilon)+\frac{\dot \rho(t)}{\rho(t)}]\Phi\leq
\frac{\beta(\rho(t)\epsilon)}{{\mathcal{R}}_0}\Phi,
\ \ \  (y,t)\in[0,L]\times[0,T].
$$

Similar as in \eqref{c08}, it is clear that
\begin{equation*}
\frac{\int_0^T\beta(0)dt}{\int_0^T\gamma(0)dt}
\leq {\mathcal{R}}_0\leq\frac{\int_0^T\beta(\rho(t)\epsilon)dt}{\int_0^T\gamma(\rho(t)\epsilon)dt}.
\end{equation*}
The result $(i)$ is given by making $\epsilon\rightarrow0$.

We now prove $(ii)$. Assume that
$$
\lim\limits_{z\to \infty}\beta(z)=\beta_{\infty}\ \ \mbox{and} \ \lim\limits_{z\to \infty}\gamma(z)=\gamma_{\infty}.
$$
Hence, the following results hold that
\begin{equation}
\lim\limits_{L\to \infty}\frac{1}{L}\int_0^L\beta(\rho(t)y)dy=\beta_{\infty}
\label{cth301}
\end{equation}
and
\begin{equation}
\lim\limits_{L\to \infty}\frac{1}{L}\int_0^L\gamma(\rho(t)y)dy=\gamma_{\infty}.
\label{cth302}
\end{equation}
In fact, $\lim\limits_{z\to \infty}\beta(z)=\beta_{\infty}$, so for any small $\epsilon>0$, there exists
a large enough $M>0$, when $y>M$
$$
|\beta(\rho(t)y)-\beta_{\infty}|<\epsilon
$$
for given $L$, there exists $L^{*}>M$, while $L>L^*$
$$
\frac{\int_0^M(\beta(\rho(t)y)-\beta_{\infty})dy}{L}<\epsilon,
$$
and
$$\begin{array}{llllll}
&&\frac{1}{L}\int_0^L(\beta(\rho(t)y)-\beta_{\infty})dy\\[2mm]
&=&\frac{1}{L}[\int_0^M(\beta(\rho(t)y)-\beta_{\infty})dy+\int_M^L(\beta(\rho(t)y)-\beta_{\infty})dy]\\[2mm]
&\leq&\frac{\int_0^M(\beta(\rho(t)y)-\beta_{\infty})dy}{L}+\frac{\epsilon(L-M)}{L}\\[2mm]
&<&2\epsilon.
\end{array}$$
Therefore, \eqref{cth301} holds. In the same way, \eqref{cth302} can be proved.

Let $s=\frac{y}{L}, 0<s<1$, then $\Phi(y,t)=\Phi(sL,t)$. Define
$$
\varphi(s,t):=\Phi(y,t), \ \ \ 0<s<1,
$$
then $\varphi(s,t)$ satisfies
\begin{eqnarray}
\left\{
\begin{array}{ll}
\varphi_{t}-\frac{d_IL^2}{\rho^2(t)} \varphi_{ss}+[\gamma(\rho(t)sL)+\frac{\dot \rho(t)}{\rho(t)}]\varphi=\frac{\beta(\rho(t)sL)}{{\mathcal{R}}_0}\varphi,\; &\ 0<s<1, 0<t\leq T, \\[2mm]
\frac{\partial \varphi}{\partial \nu}=0,\ & \ s=0,1, 0<t\leq T,\\[2mm]
\varphi(s,0)=\varphi(s,T),\; &\ 0\leq s\leq 1.
\end{array} \right.
\label{cth303}
\end{eqnarray}
We now normalize $\varphi(s,t)$ such that
\begin{equation}
\int_0^T\int_0^1\varphi^{2}(s,t)dsdt=1.
\label{cth304}
\end{equation}
We multiply the equation of \eqref{cth303} by $\varphi(s,t)$ and integrate over $(0,1)\times(0,T)$ to obtain
$$
\int_0^T\int_0^1\frac{d_IL^{2}}{\rho^{2}(t)}{\varphi^{2}_{s}}dsdt=\int_0^T\int_0^1[\frac{\beta(\rho(t)sL)}
{{\mathcal{R}}_0}-\gamma(\rho(t)sL)-\frac{\dot{\rho(t)}}{\rho(t)}]\varphi^{2}dsdt.
$$
Combining  \eqref{cth304} with \eqref{c08} yields
$$
\frac{d_IL^{2}}{(\rho^{M})^2}\int_0^{T}\int_0^{1}\varphi^{2}_{s}dsdt\leq\int_0^{T}\int_0^{1}\frac{d_IL^2}
{\rho^{2}(t)}\varphi^{2}_{s}dsdt\leq C,
$$
i.e.,
\begin{equation}
\int_0^{T}\int_0^{1}\varphi^{2}_{s}dsdt\leq\frac{C(\rho^{M})^2}{d_IL^2},
\label{cth305}
\end{equation}
where $C$ is positive constant and independent of $L$.

On the other hand, define
$$
\varphi_{*}(t):=\int_0^{1}\varphi(s,t)ds \ \ \mbox{and} \ \ \varphi^{*}(s,t):=\varphi(s,t)-\varphi_{*}(t).
$$
It is easy to see that
$$
\int_0^{1}\varphi^{*}(s,t)ds=0 \ \ \mbox{for all} \ \ t\in [0,T].
$$
Applying for the well-known Poincar$\acute{e}$ inequality, it follows that
$$
\int_0^{1}(\varphi^{*}(s,t))^{2}ds \leq C\int_0^{1}(\varphi^{*}_{s})^{2}ds \ \ \mbox{for all} \ \ t\in [0,T].
$$
Making use of \eqref{cth305} and $\varphi^{*}_{s}=\varphi_{s}$, it gives
\begin{equation}
\int_0^{T}\int_0^{1}(\varphi^{*})^{2}dsdt\leq\frac{C(\rho^{M})^2}{d_IL^2}.
\label{cth306}
\end{equation}
Using the H$\ddot{o}$lder inequality for \eqref{cth305} and \eqref{cth306} gives
\begin{equation}
\int_0^{T}\int_0^{1}|\varphi^{*}|dsdt\leq C\rho^{M}(\frac{1}{d_IL^2})^{\frac{1}{2}},\ \ \
\int_0^{T}\int_0^{1}|\varphi^{*}_{s}|dsdt\leq C\rho^{M}(\frac{1}{d_IL^2})^{\frac{1}{2}}.
\label{cth307}
\end{equation}
In addition, integrating the equation of \eqref{cth303} over $(0,1)$, we have
{\small \begin{equation}
\frac{d\varphi_{*}}{dt}=\varphi_{*}\int_0^{1}[\frac{\beta(\rho(t)sL)}{{\mathcal{R}}_0}-\gamma(\rho(t)sL)
-\frac{\dot \rho(t)}{\rho(t)}]ds+\int_0^{1}[\frac{\beta(\rho(t)sL)}{{\mathcal{R}}_0}-\gamma(\rho(t)sL)-
\frac{\dot \rho(t)}{\rho(t)}]\varphi^*ds.
\label{cth308}
\end{equation}}
Solving the ODE \eqref{cth308} yields
\begin{equation}
\varphi_{*}(t)=\varphi_{*}(0)\exp(\int_0^t\int_0^1(\frac{\beta(\rho(t)sL)}{{\mathcal{R}}_0}-\gamma(\rho(t)sL)
-\frac{\dot \rho(t)}{\rho(t)})dsdt)+A(t),
\label{cth309}
\end{equation}
where
{\small $$
A(t)=\int_0^t\int_0^1\exp(\int_\tau^t\int_0^1(\frac{\beta(\rho(t)sL)}{{\mathcal{R}}_0}-\gamma(\rho(t)sL)-\frac{\dot \rho(t)}{\rho(t)})dsdt)[\frac{\beta(\rho(t)sL)}{{\mathcal{R}}_0}-\gamma(\rho(t)sL)-\frac{\dot \rho(t)}{\rho(t)}]\varphi^*dsdt
$$}
for $\tau, t\in(0,T)$.
According to \eqref{c08}, there exists a constant $M>0$ such that
$$
\exp(\int_\tau^t\int_0^1(\frac{\beta(\rho(t)sL)}{{\mathcal{R}}_0}-\gamma(\rho(t)sL)-\frac{\dot \rho(t)}{\rho(t)})dsdt\leq M.
$$
In view of \eqref{cth307} and above inequality, we know
{\small $$
A(t)\leq M\int_0^T\int_0^1|[\frac{\beta(\rho(t)sL)}{{\mathcal{R}}_0}-\gamma(\rho(t)sL)-\frac{\dot \rho(t)}{\rho(t)}]
\varphi^*|dsdt\rightarrow0 \ \mbox{uniformly on}\, [0,T] \, \mbox{as} \, L\rightarrow\infty.
$$}
Making $t=T$ in \eqref{cth309} and using $\varphi_{*}(0)=\varphi_{*}(T)$ yield
\begin{equation}
\varphi_{*}(0)[1-\exp(\int_0^T\int_0^1(\frac{\beta(\rho(t)sL)}{{\mathcal{R}}_0}-\gamma(\rho(t)sL))dsdt)]=A(T).
\label{cth310}
\end{equation}
Recalling that $A(T)\rightarrow 0$ as $L\rightarrow\infty$, we next conclude that $\varphi_{*}(0)\nrightarrow0$ as $L\rightarrow\infty$.

In fact, if $\varphi_{*}(0)\rightarrow0$ as $L\rightarrow\infty$. In light of \eqref{cth309}, we easily see that
\begin{equation}
\varphi_{*}(t)\rightarrow0 \ \ \mbox{uniformly on}\ \  [0,T] \ \ \mbox{as} \ \ L\rightarrow\infty.
\label{cth311}
\end{equation}
Combining \eqref{cth307} with \eqref{cth311} ensures
$$
\int_0^T\int_0^1\varphi dsdt\rightarrow 0 \ \ \mbox{as} \ \ \ L\rightarrow \infty,
$$
which is a contradiction with \eqref{cth304}.

Now substituting $y=sL$ into \eqref{cth310}, then it becomes
$$
\varphi_{*}(0)[1-\exp(\int_0^T\int_0^L\frac{1}{L}(\frac{\beta(\rho(t)y)}{{\mathcal{R}}_0}-\gamma(\rho(t)y))dydt)]=A(T),
$$
which gives
\begin{equation*}
{\mathcal{R}}_0=\frac{\frac{1}{L}\int_0^T\int_0^L\beta(\rho(t)y)dydt}{\frac{1}{L}\int_0^T\int_0^L\gamma(\rho(t)y)dydt
+\ln(1-\frac{A(T)}{\varphi_{*}(0)})}.
\label{cth312}
\end{equation*}
Due to $A(T)\rightarrow 0$ and $\varphi_{*}(0)\nrightarrow0$ as $L\rightarrow\infty$, then
$$
\ln(1-\frac{A(T)}{\varphi_{*}(0)})\rightarrow 0 \ \ \ \ \mbox{as} \ \ \ L\rightarrow\infty.
$$
In light of above result, \eqref{cth301} and \eqref{cth302}, it is obvious that
$$
{\mathcal{R}}_0=\frac{\beta_{\infty}}{\gamma_{\infty}} \ \ \ \ \mbox{as} \ \ \ L\rightarrow\infty.
$$

\epf

\begin{rmk}
\label{ma}
Assume that $\gamma_{y}(\rho(t)y)\geq0$ and $\beta_{y}(\rho(t)y)<0$ for any $t>0$. From the above proof, we can find that $\Phi_{y}(y,t)<0$,
Theorem \ref{main3} still holds, and $(i)$ of Theorem \ref{main1} and $(ii)$ of Theorem \ref{main2} hold, but
${\mathcal{R}}_0$ is strictly monotone decreasing in $L$ for any given $d_I>0$, and $\lim\limits_{d_I\to 0}{{\mathcal{R}}_0}=\frac{\beta(0)}{\gamma(0)}$.
\end{rmk}

In some special cases, the explicit expression of ${\mathcal{R}}_0$ can be obtained. In fact, suppose that the coefficients $\beta(\rho(t)y)\equiv\widehat{\beta}$, $\gamma(\rho(t)y)=\frac{c(y)}{\rho^2(t)}+g(t)$,
$a(\rho(t)y)\equiv \widehat{a}$ and $b(\rho(t)y)\equiv\widehat{b}$ in problem \eqref{a07}, where $\widehat{\beta}$,
$\widehat{a}$ and $\widehat{b}$ are positive constants and $g(t)$ is T-periodic in time. Set $\lambda^*$
is the principal eigenvalue of the following problem
\begin{eqnarray}
\left\{
\begin{array}{ll}
-d_{I}\Delta\varpi+c(y)\varpi=\lambda\varpi,\; &\ y\in \Omega(0), \\[2mm]
\frac{\partial \varpi}{\partial \nu}=0,\ & \ y\in \partial \Omega(0).
\end{array} \right.
\label{d03}
\end{eqnarray}
In view of \eqref{c01}, let $\omega=e^{\int_{0}^{t}g(s)ds}\phi$, then $\omega$ satisfies
\begin{eqnarray}
\left\{
\begin{array}{ll}
\omega_{t}-\frac{d_I}{\rho^2(t)} \Delta \omega+\frac{c(y)}{\rho^2(t)}\omega=\frac{\widehat{\beta}}{R_0}\omega-
\frac{n\dot \rho(t)}{\rho(t)}\omega,\; &\ y\in \Omega(0), t>0, \\[2mm]
\frac{\partial \omega}{\partial \nu}=0,\ & \ y\in \partial \Omega(0), t>0,\\[2mm]
\omega(y,0)=\omega(y,T),\; &\ y\in \overline \Omega(0).
\end{array} \right.
\label{d04}
\end{eqnarray}
Straightforward calculation alleges
\begin{eqnarray}
{\mathcal{R}}_0=\frac{\int_0^{T}\widehat{\beta} dt}{\int_0^{T}\frac{\lambda^{*}}{\rho^2(t)}dt},
\label{d05}
\end{eqnarray}
which implies that ${\mathcal{R}}_0$ decreases with respect
to $\overline {\rho^{-2}}$ (:=$\frac 1T\int^T_0 \frac 1{\rho^2(t)}dt$).

\section{\bf Simulation and discussion}

\begin{figure}[ht]
\centering
\subfigure[]{ {
\includegraphics[width=0.28\textwidth]{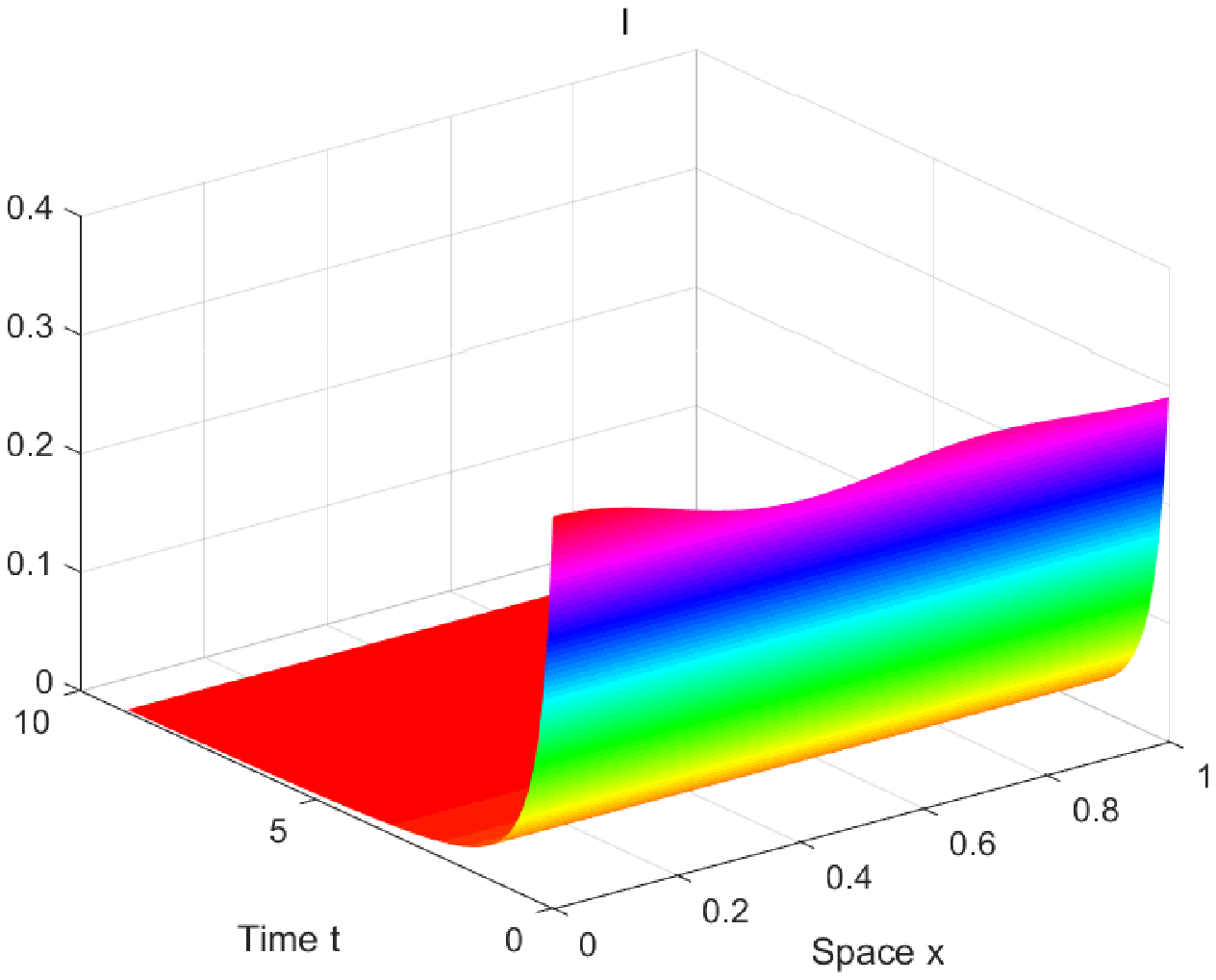}
} }
\subfigure[]{ {
\includegraphics[width=0.28\textwidth]{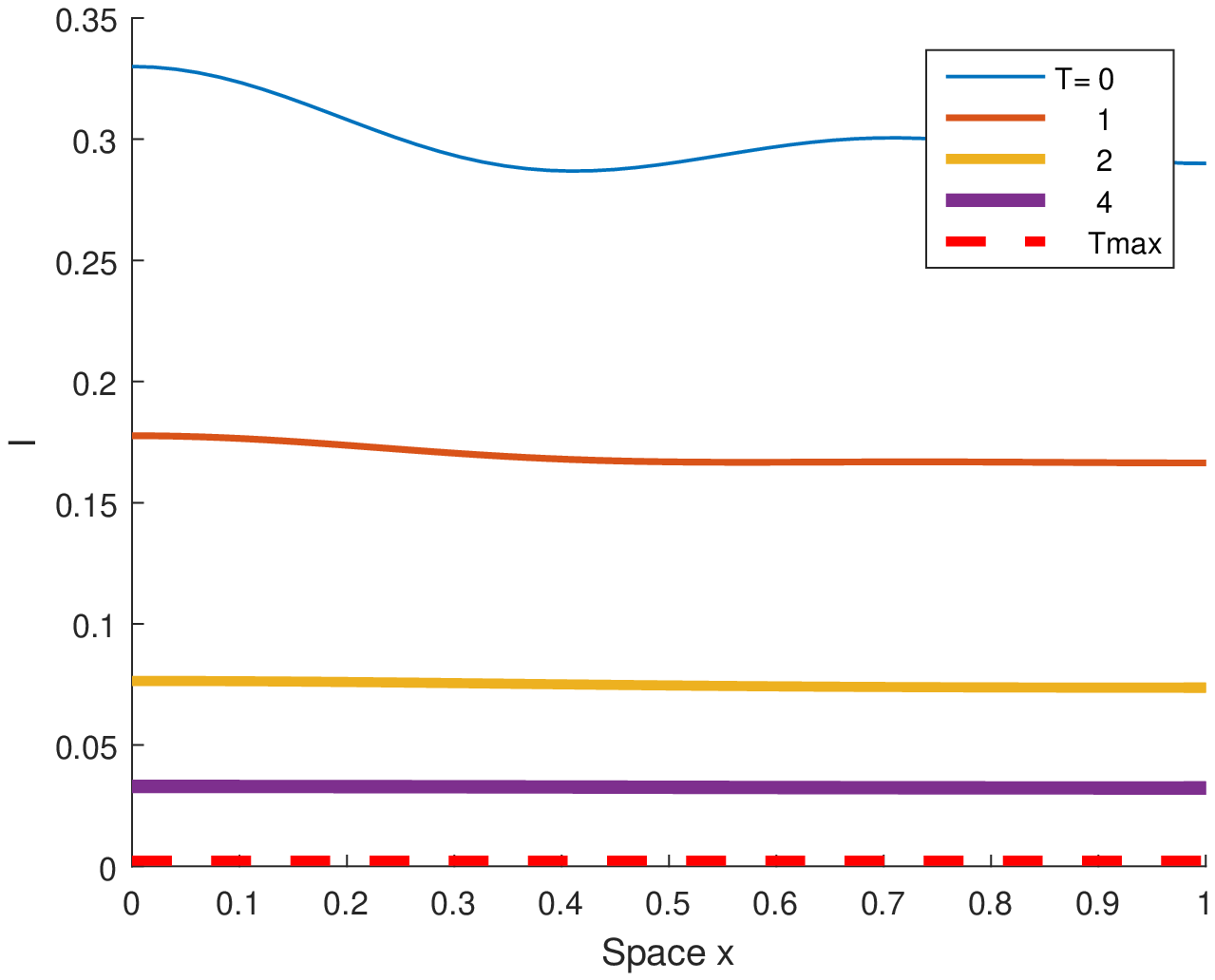}
} }
\subfigure[]{ {
\includegraphics[width=0.28\textwidth]{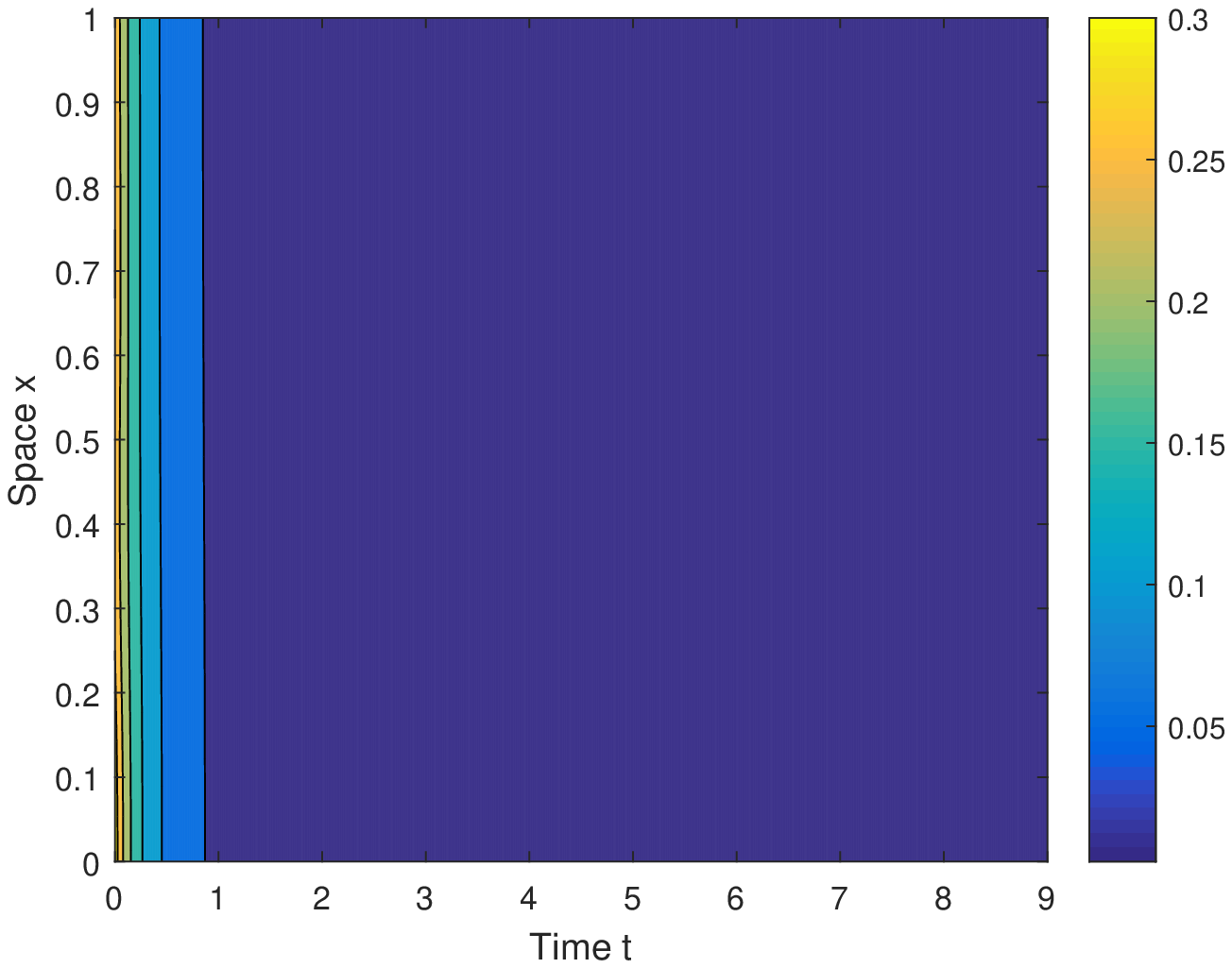}
} }
\caption{\scriptsize $c=6.96$ and $\rho_1(t)\equiv 1$. The domain is fixed and ${\mathcal{R}}_0<1$.
Graph $(a)$ shows that infected individual $I$ decays to $0$. Graphs $(b)$ and $(c)$ are the cross-sectional
view and contour map respectively.}
\label{tu1}
\end{figure}

\begin{figure}[ht]
\centering
\subfigure[]{ {
\includegraphics[width=0.28\textwidth]{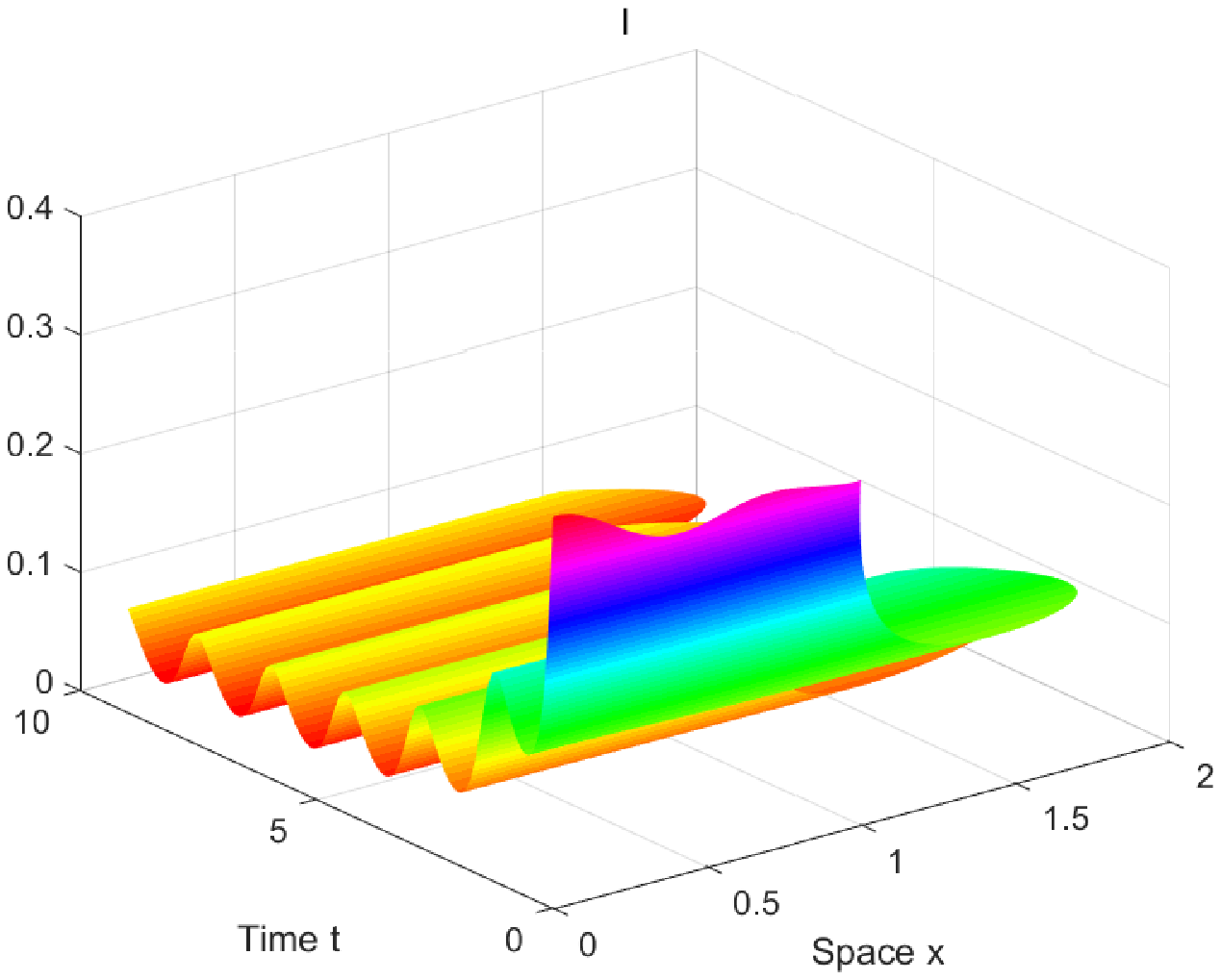}
} }
\subfigure[]{ {
\includegraphics[width=0.28\textwidth]{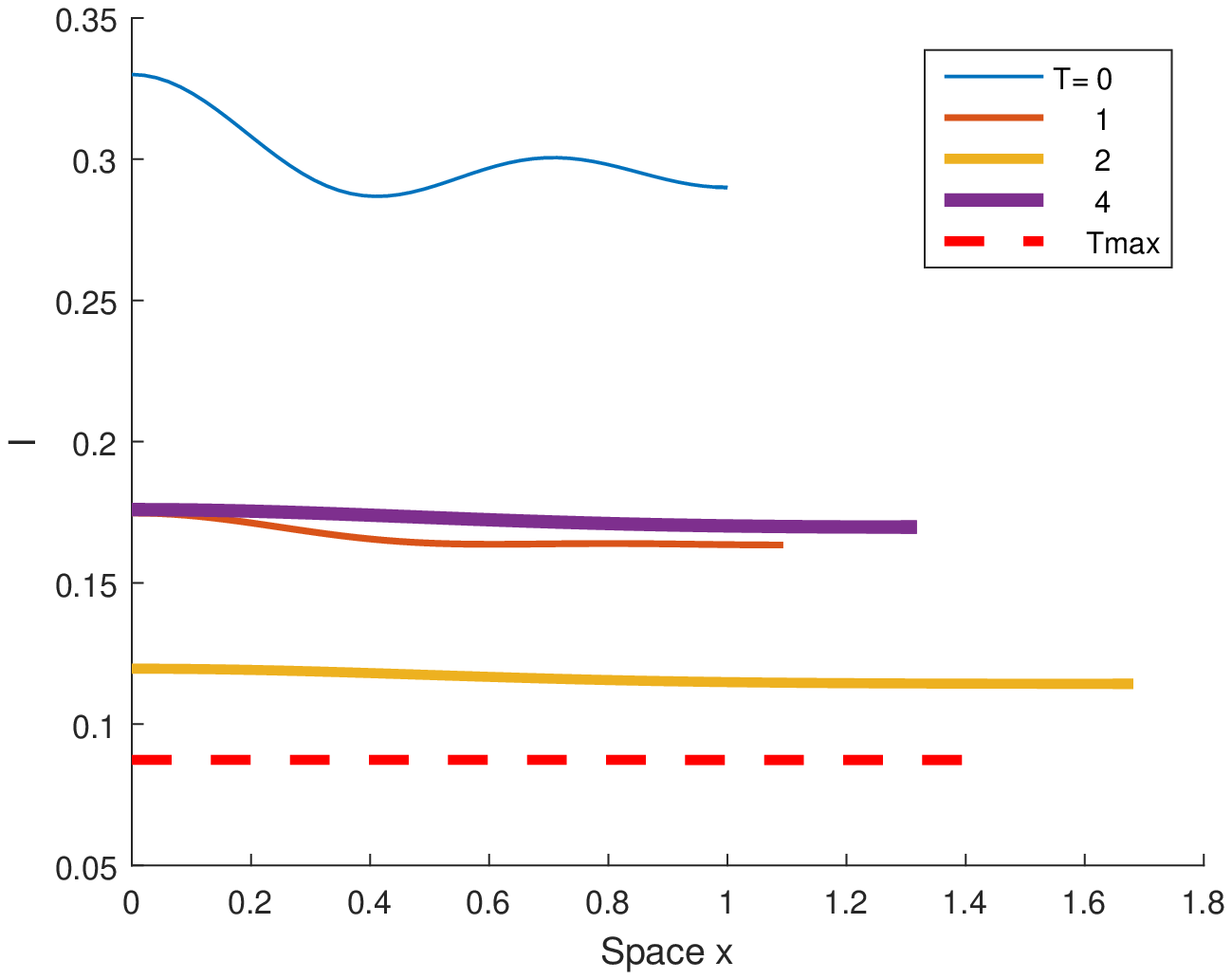}
} }
\subfigure[]{ {
\includegraphics[width=0.28\textwidth]{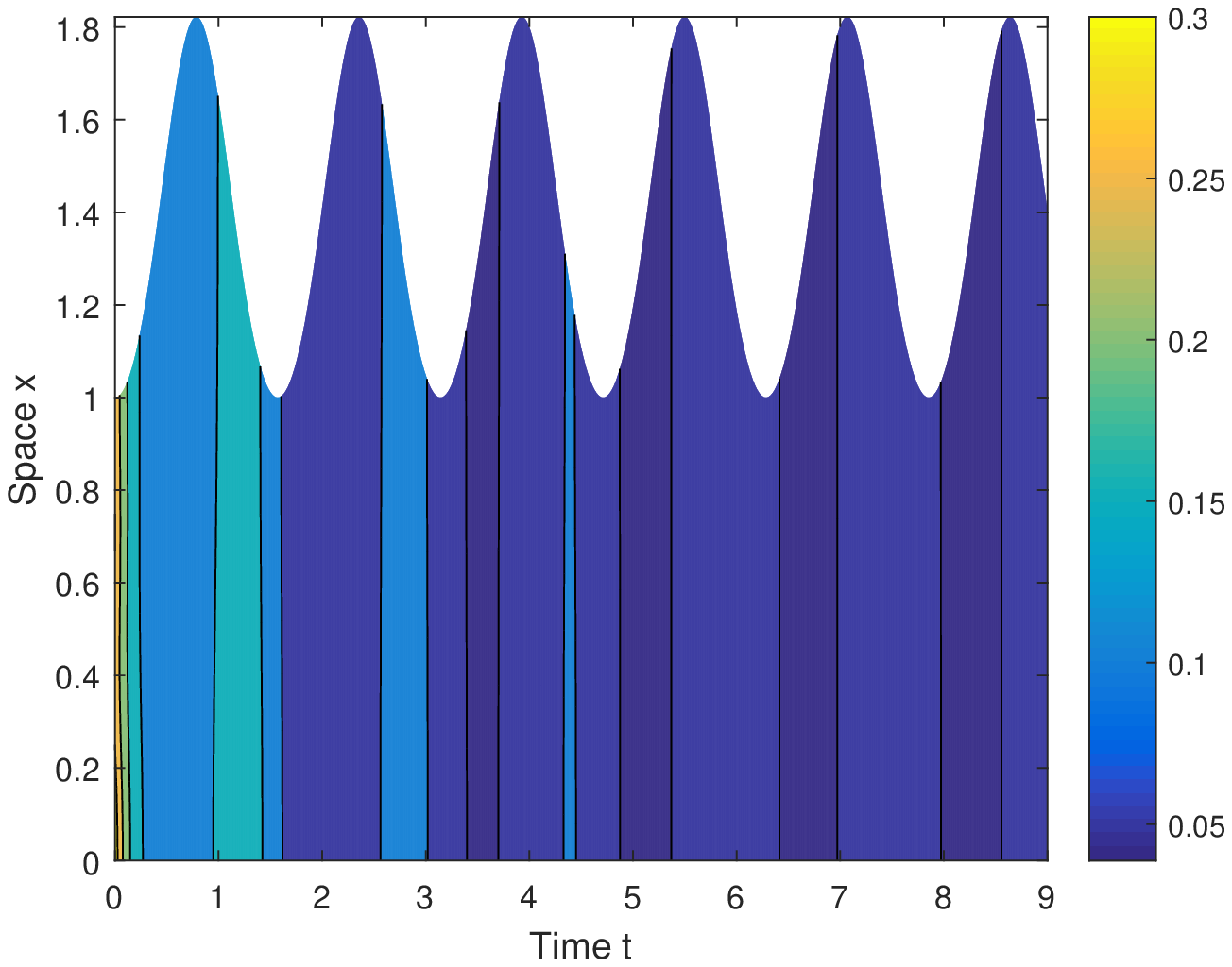}
} }
\caption{\scriptsize $c=6.96$ and $\rho_2(t)=e^{0.3(1-\cos(4t))}$. The domain is evolving with a larger
evolution rate $\rho_2(t)$ and ${\mathcal{R}}_0>1$. Graph $(a)$ shows that infected individual $I$
stabilizes to a positive periodic steady state. Graphs $(b)$ and $(c)$, which are the cross-sectional
view and contour map respectively, present the periodic evolution of the domain.}
\label{tu2}
\end{figure}

In this section, we carry out numerical simulations for problem \eqref{a07}-\eqref{a08} to explain
the theoretical results. 
Suppose that
\begin{equation*}
\left.
\begin{array}{lll}
d_S=0.01,\ d_I=0.1, \  a=1,\  b=10, \ \beta=7, \  \Omega(0)=(0, 1), \ \gamma=\frac{c}{\rho^{2}(t)}, \\[7pt]
S_0(y) = 0.3+0.01 \cos(\pi y)+0.01\cos(4\pi y),\\
I_0(y) = 0.3+0.01\cos(\pi y)+0.01\cos(2\pi y)+0.01\cos(3\pi y)
\end{array}
\right.
\end{equation*}
in \eqref{a07}-\eqref{a08} and subsequently $\lambda^*=c+d_I\pi^{2}$ in \eqref{d03}, then the asymptotic
behaviors of the solution to problem \eqref{a07}-\eqref{a08} are shown by choosing different $\rho(t)$ and $c$.
\begin{exm}
Taking $c=6.96$. We first consider that the habitat is fixed, i.e., $\rho_1(t)\equiv 1$. Calculations show that
$$
{\mathcal{R}}_0(\rho_1)=\frac{\beta}{\frac{c+d_I\pi^{2}}{\rho^2_1}}=
\frac{7}{6.69+0.1\times\pi^{2}}
\approx0.8808<1.$$
It is easy to see from Fig. \ref{tu1} that the infected individual $I$ decays to zero.

We now choose $\rho_2(t)=e^{0.3(1-\cos(4t))}$, it follows from \eqref{d05} that
$$\overline{\rho^{-2}_2}=\frac{2}{\pi}\int_0^{\frac{\pi}{2}}e^{1.2(\cos(4t)-1)}dt\approx 0.5593$$
and
$${\mathcal{R}}_0(\rho_2)=\frac{\int_0^{\frac{\pi}{2}}\beta dt}{\int_0^{\frac{\pi}{2}}\frac{c+d_I\pi^{2}}{\rho^2_2(t)}dt}=
\frac{\beta}{(c+d_I\pi^{2})\overline{\rho^{-2}_2}}
\approx 1.5749>1.$$
It is easy to see from Fig. \ref{tu2} that $I$ stabilizes to a positive periodic steady state.

This example shows that the infected individual vanishes on a fixed domain, but persists on a periodically
evolving domain.
\end{exm}


\begin{figure}[ht]
\centering
\subfigure[]{ {
\includegraphics[width=0.28\textwidth]{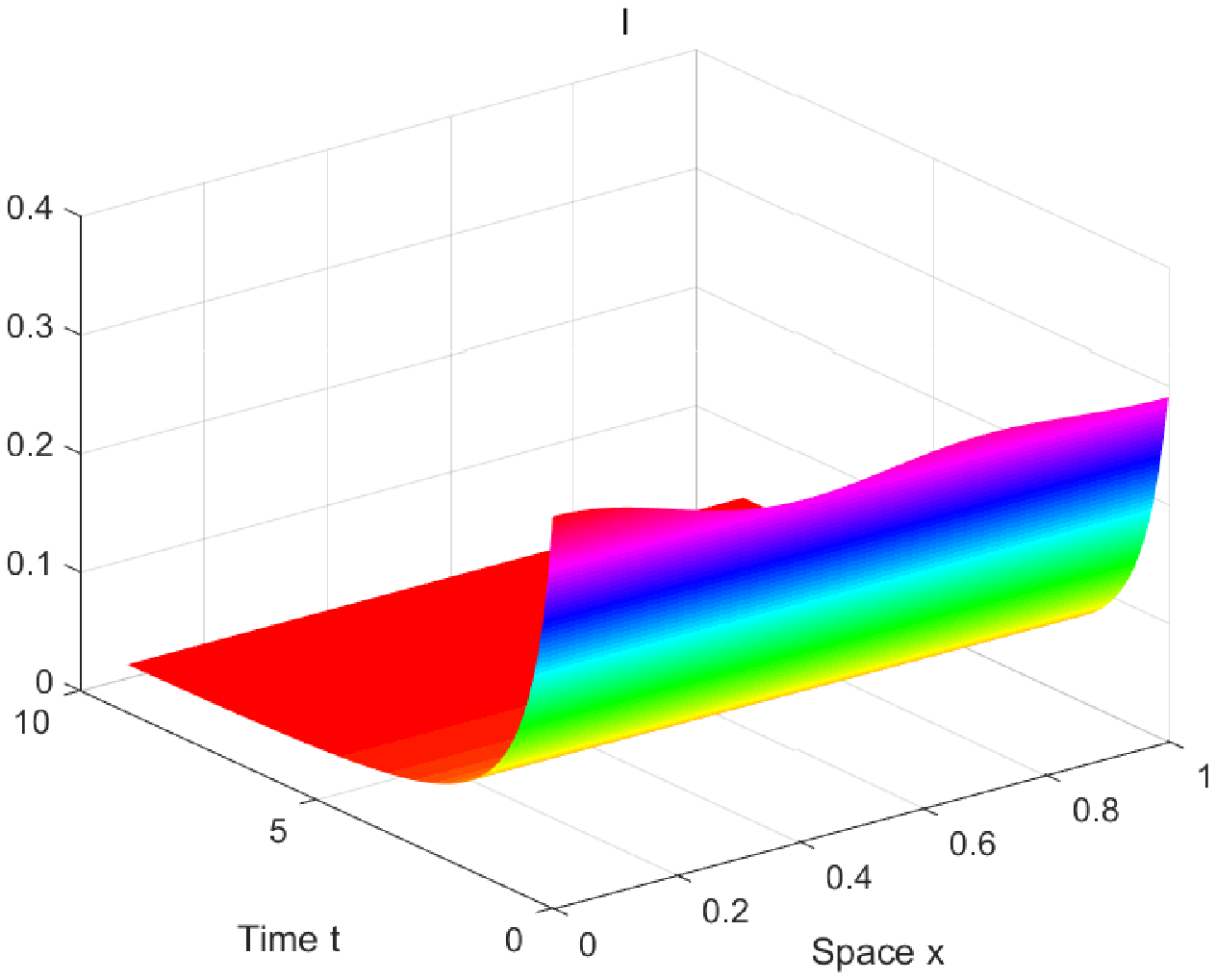}
} }
\subfigure[]{ {
\includegraphics[width=0.28\textwidth]{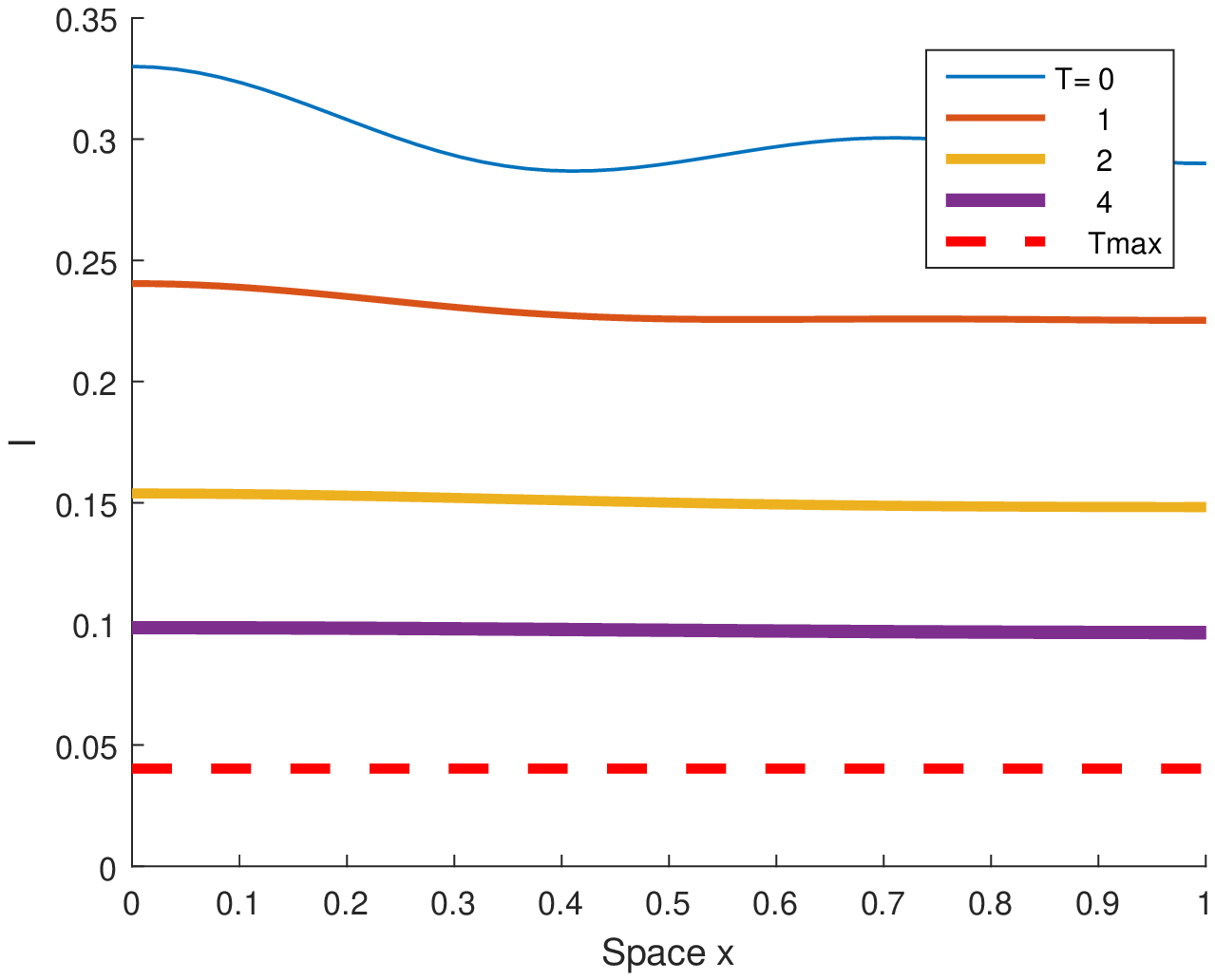}
} }
\subfigure[]{ {
\includegraphics[width=0.28\textwidth]{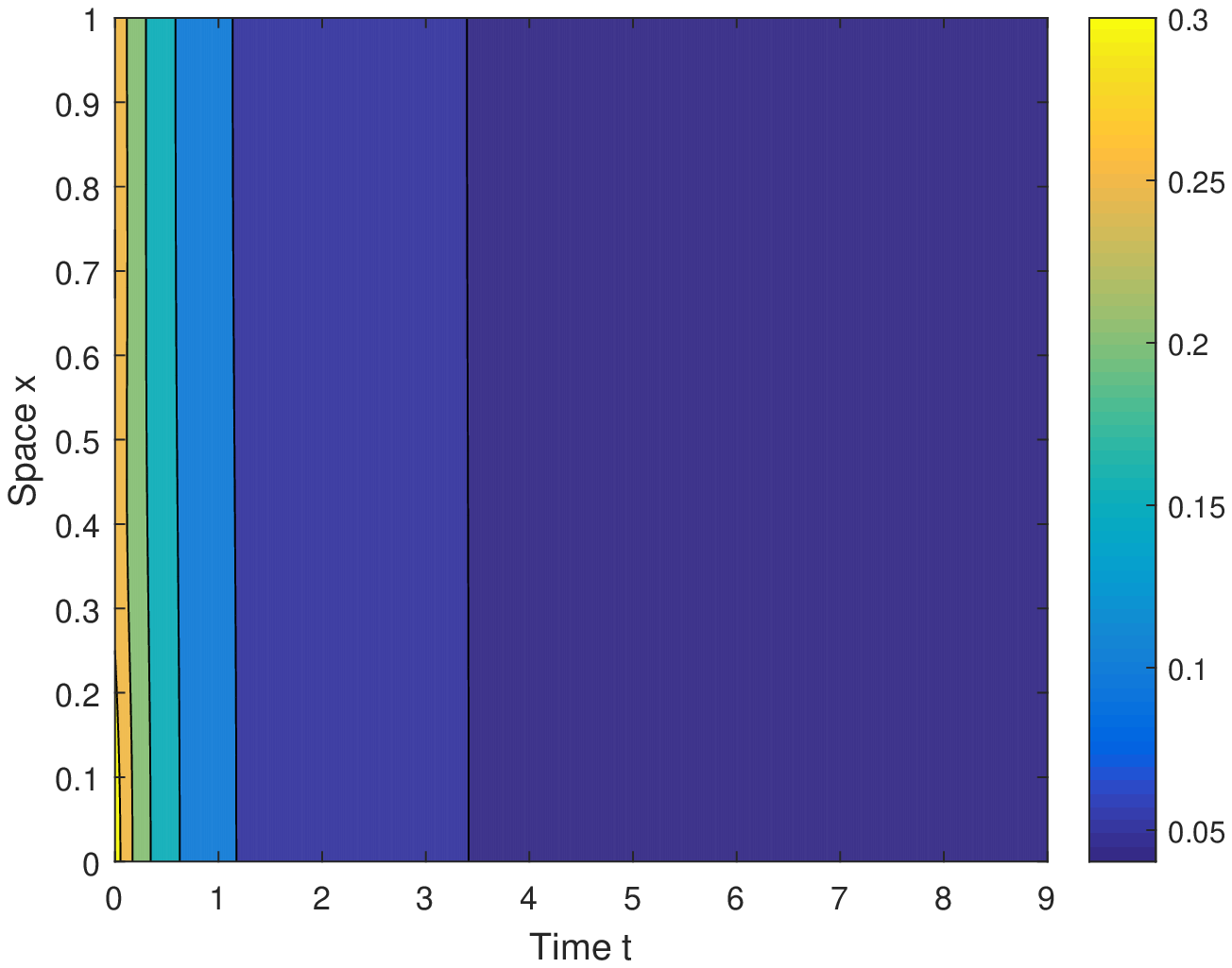}
} }
\caption{\scriptsize $c=5$ and $\rho_3(t)\equiv 1$. The domain is fixed and ${\mathcal{R}}_0>1$. Graph $(a)$ shows that
infected individual $I$ tends to a positive periodic steady state. Graphs $(b)$ and $(c)$ are the cross-sectional view
and contour map respectively. }
\label{tu3}
\end{figure}

\begin{figure}[ht]
\centering
\subfigure[]{ {
\includegraphics[width=0.28\textwidth]{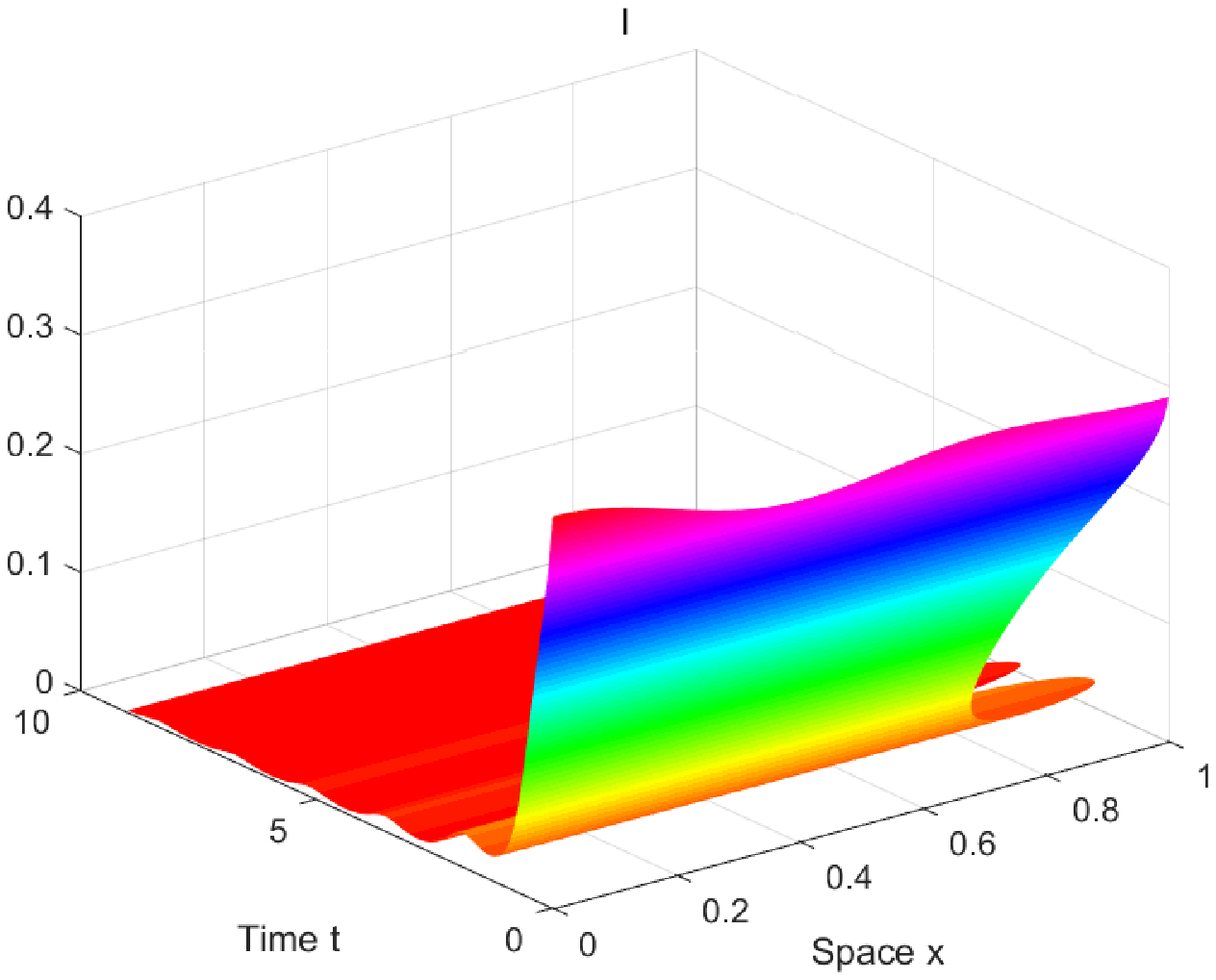}
} }
\subfigure[]{ {
\includegraphics[width=0.28\textwidth]{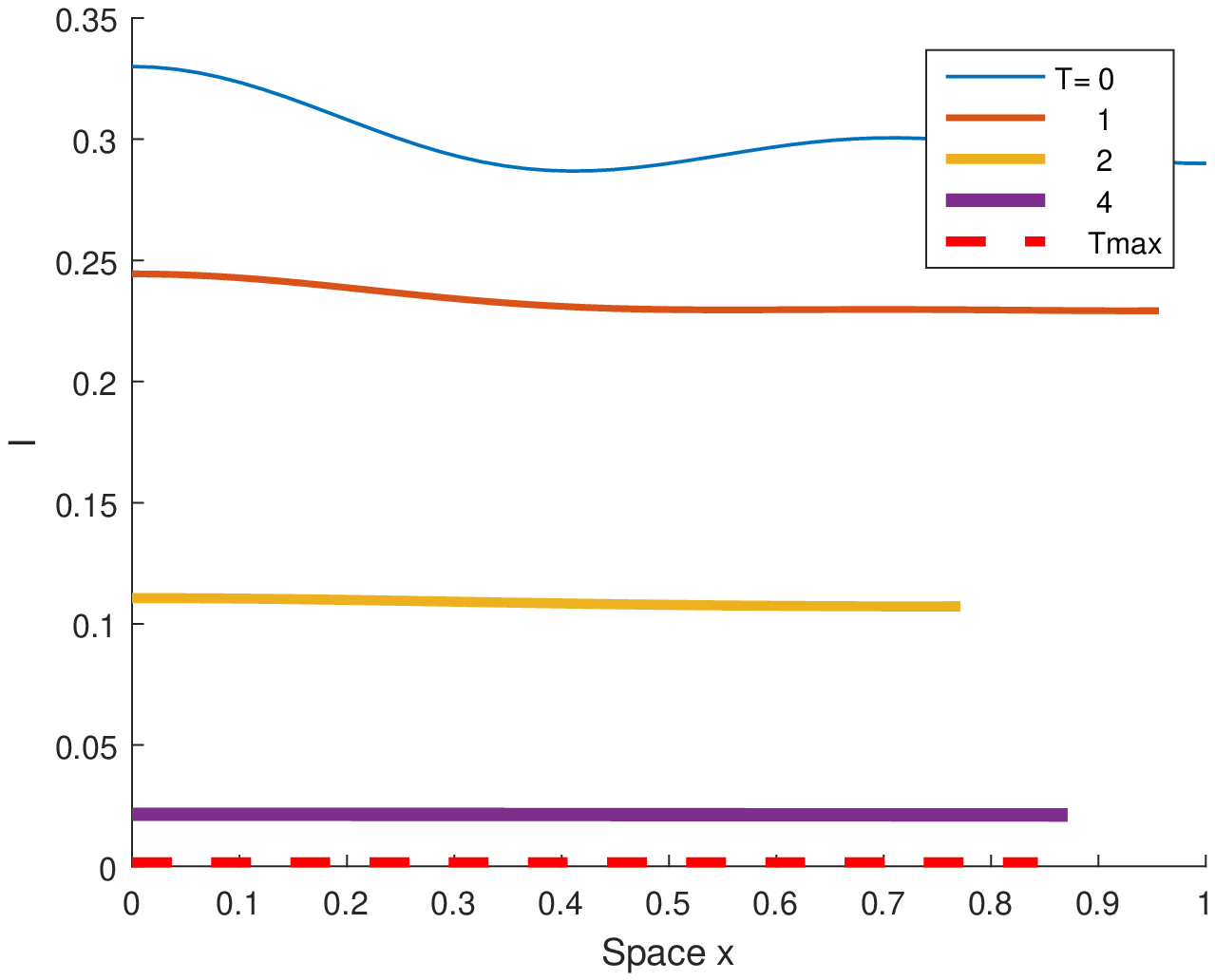}
} }
\subfigure[]{ {
\includegraphics[width=0.28\textwidth]{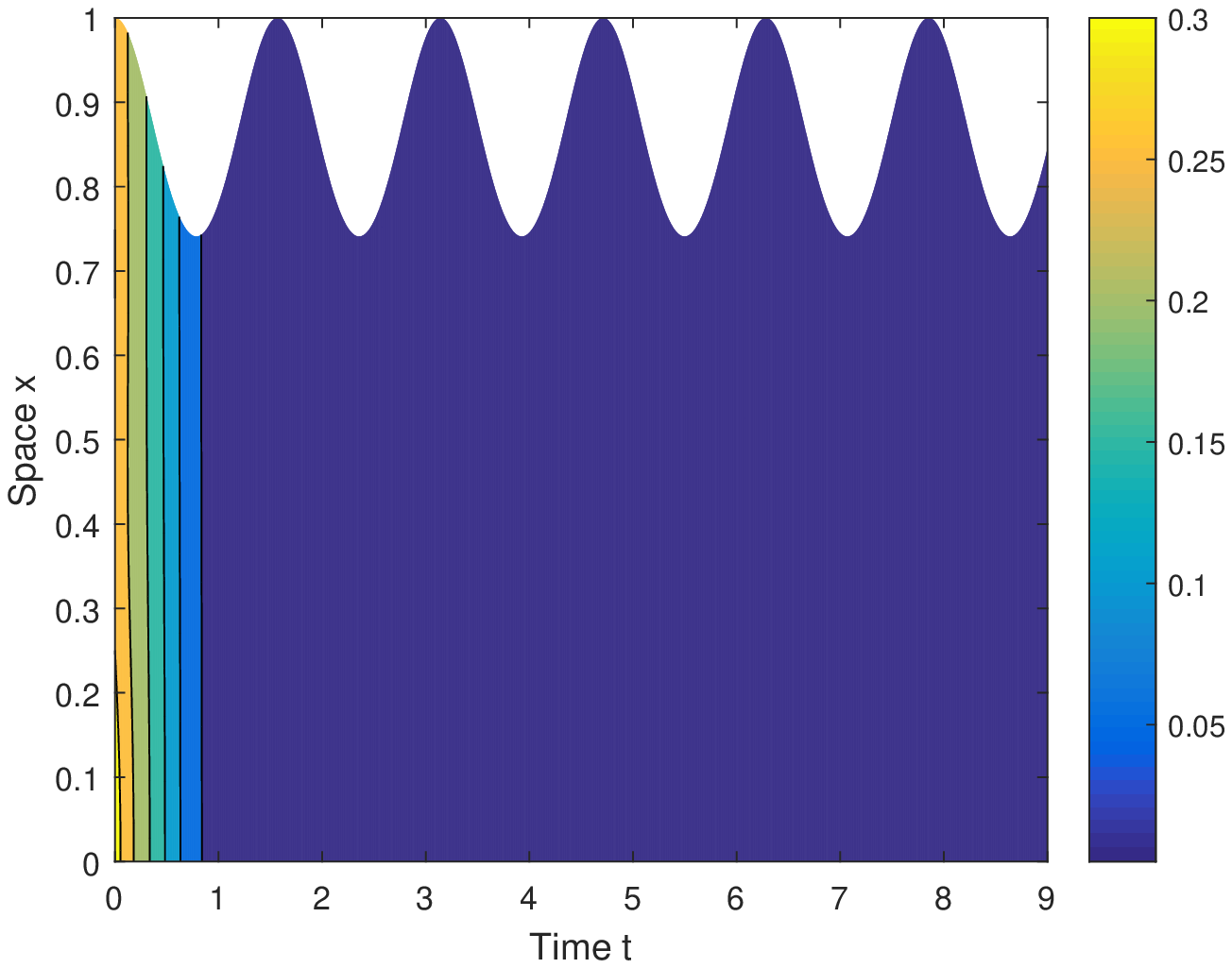}
} }
\caption{\scriptsize $c=5$ and $\rho_4(t)=e^{-0.15(1-\cos(4t))}$. The domain is evolving with a smaller
evolution rate $\rho_4(t)$ and ${\mathcal{R}}_0<1$. Graph $(a)$ shows that infected individual $I$ decays
to zero. Graphs $(b)$ and $(c)$, which are the cross-sectional view and contour map respectively, present
the periodic evolution of the domain.
}
\label{tu4}
\end{figure}

\begin{exm}
Letting $c=5$. At first, we choose $\rho_3(t)\equiv 1$, which implies that the habitat is a fixed domain.
Straightforward calculation alleges
$$
{\mathcal{R}}_0(\rho_3)=\frac{\beta}{\frac{c+d_I\pi^{2}}{\rho^2_1}}=
\frac{7}{5+0.1\times\pi^{2}}
\approx1.1692>1.
$$
It is easy to see from Fig. \ref{tu3} that $I$ stabilizes to a positive periodic steady state.

Next, we choose $\rho_4(t)=e^{-0.15(1-\cos(4t))}$, it follows from \eqref{d05} that
$$\overline{\rho^{-2}_4}=\frac{2}{\pi}\int_0^{\frac{\pi}{2}}e^{0.6(1-\cos(4t))}dt\approx 1.3804$$
and
$${\mathcal{R}}_0(\rho_4)=\frac{\int_0^{\frac{\pi}{2}}\beta dt}{\int_0^{\frac{\pi}{2}}
\frac{c+d_I\pi^{2}}{\rho^2_2(t)}dt}=\frac{\beta}{(c+d_I\pi^{2})\overline{\rho^{-2}_2}}
\approx 0.8470<1.$$
It is easy to see from Fig. \ref{tu4} that $I$ decays to zero and the infected individual vanishes eventually.

The above example tells us that the infected individual spreads on a fixed domain, but disappears on a periodically
evolving domain.
\end{exm}

Next, we verify the monotonicity of ${\mathcal{R}}_0$ with respect to $d_I$. We fix the $\rho_{5}(t)
=e^{0.3(1-\cos(4t))}$.

\begin{figure}[ht]
\centering
\subfigure[]{ {
\includegraphics[width=0.28\textwidth]{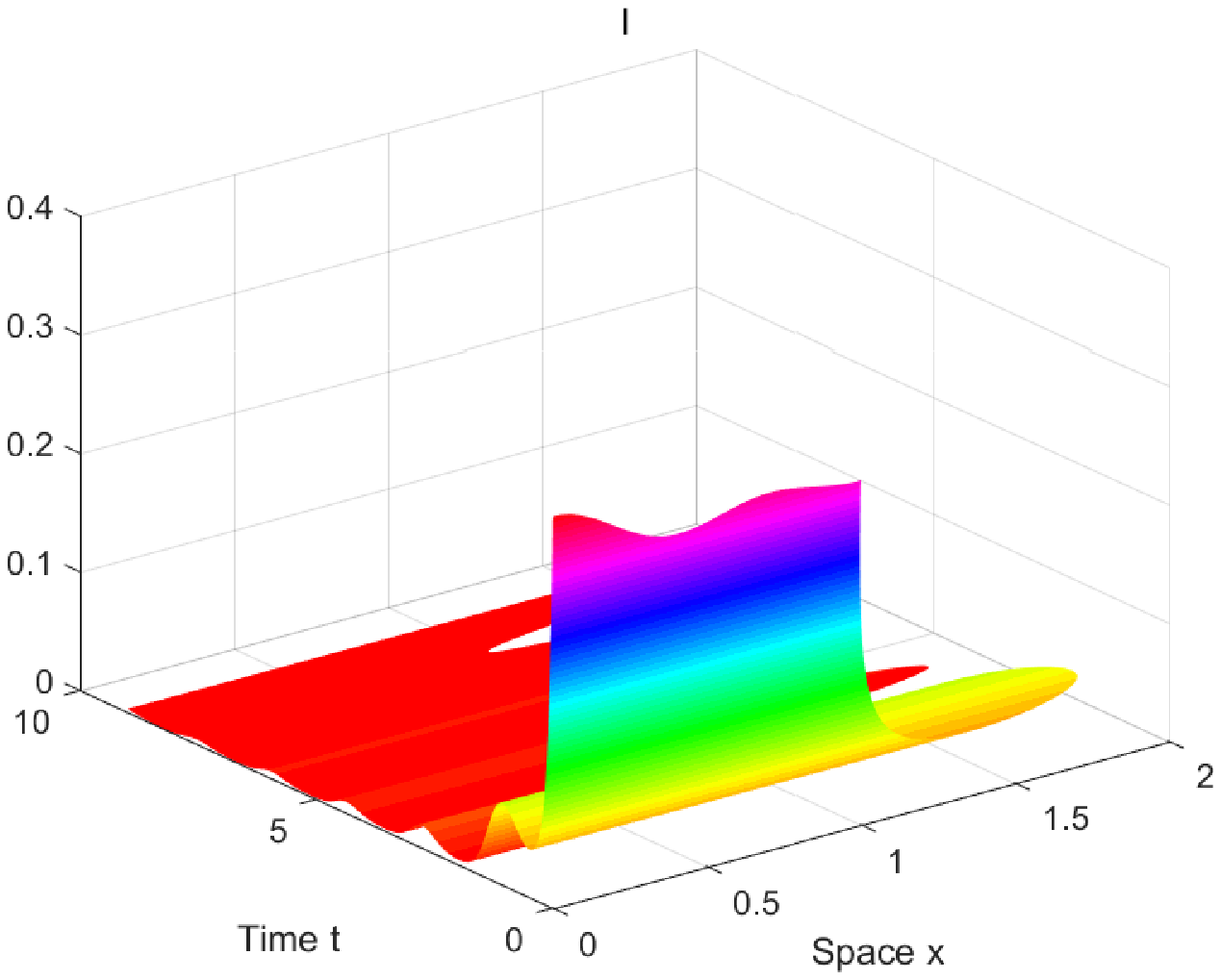}
} }
\subfigure[]{ {
\includegraphics[width=0.28\textwidth]{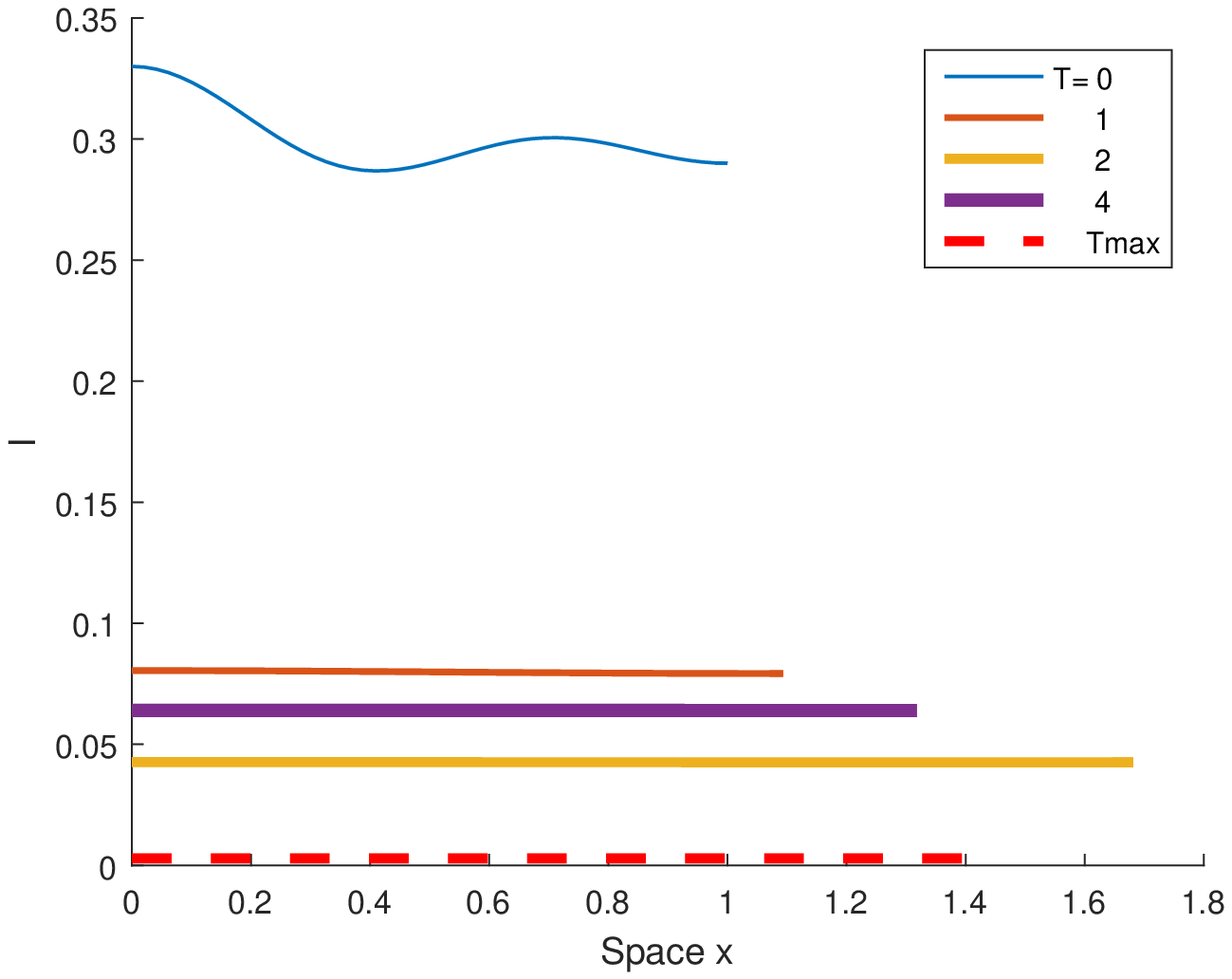}
} }
\subfigure[]{ {
\includegraphics[width=0.28\textwidth]{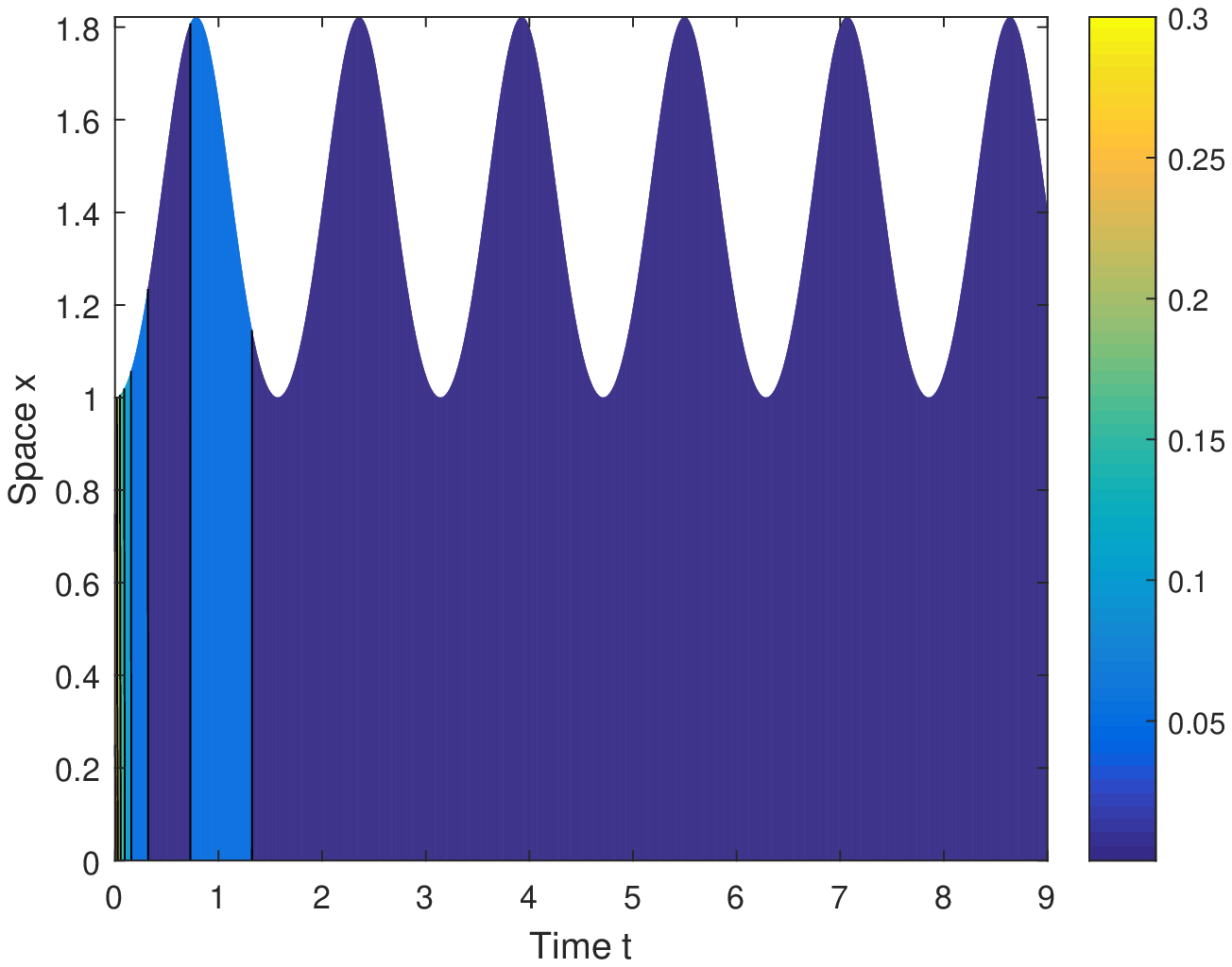}
} }
\caption{\scriptsize $c=11.8$ and $d_I=0.8$, ${\mathcal{R}}_0<1$. Graph $(a)$ shows that infected individual
$I$ decays to zero. Graphs $(b)$ and $(c)$, which are the cross-sectional view and contour map respectively.
}
\label{tu5}
\end{figure}

\begin{exm}
Letting $c=6.96$ and $d_I=0.1$. Due to \eqref{d05}, straightforward calculation alleges

$$
{\mathcal{R}}_0(\rho_5)=\frac{\int_0^{\frac{\pi}{2}}\beta dt}{\int_0^{\frac{\pi}{2}}\frac{c+d_I\pi^{2}}
{\rho^2_5(t)}dt}=
\frac{\beta}{(c+d_I\pi^{2})\overline{\rho^{-2}_5}}
\approx 1.5749>1.
$$
It is easy to see from Fig. \ref{tu2} that $I$ stabilizes to a positive periodic steady state.

We now choose $c=11.8$ and $d_I=0.8$. It follows from \eqref{d05} that

$$
{\mathcal{R}}_0(\rho_5)=\frac{\int_0^{\frac{\pi}{2}}\beta dt}{\int_0^{\frac{\pi}{2}}
\frac{c+d_I\pi^{2}}{\rho^2_5(t)}dt}=\frac{\beta}{(c+d_I\pi^{2})\overline{\rho^{-2}_5}}
\approx 0.6355<1.
$$
It is easy to see from Fig. \ref{tu5} that $I$ decays to zero and the infected individual vanishes eventually.
\end{exm}

Finally, we examine the monotonicity of ${\mathcal{R}}_0$ with respect to $L$. We redefine some coefficients
in \eqref{a07} and \eqref{a08}. Suppose that
\begin{equation*}
\left.
\begin{array}{lll}
d_S=0.1,\ d_I=0.1, \  a=1,\  b=1,  \  \rho_{6}(t)=e^{-0.2(1-\cos(2t))},\\[7pt]
S_0(y) = 0.3+0.01 \cos(\pi y)+0.01\cos(4\pi y),\\
I_0(y) = 0.3+0.01\cos(\pi y)+0.01\cos(2\pi y)+0.01\cos(3\pi y),\\
\beta=0.33+0.01e^{-0.01\rho_{6}(t)y}.
\end{array}
\right.
\end{equation*}

\begin{figure}[ht]
\centering
\subfigure[]{ {
\includegraphics[width=0.28\textwidth]{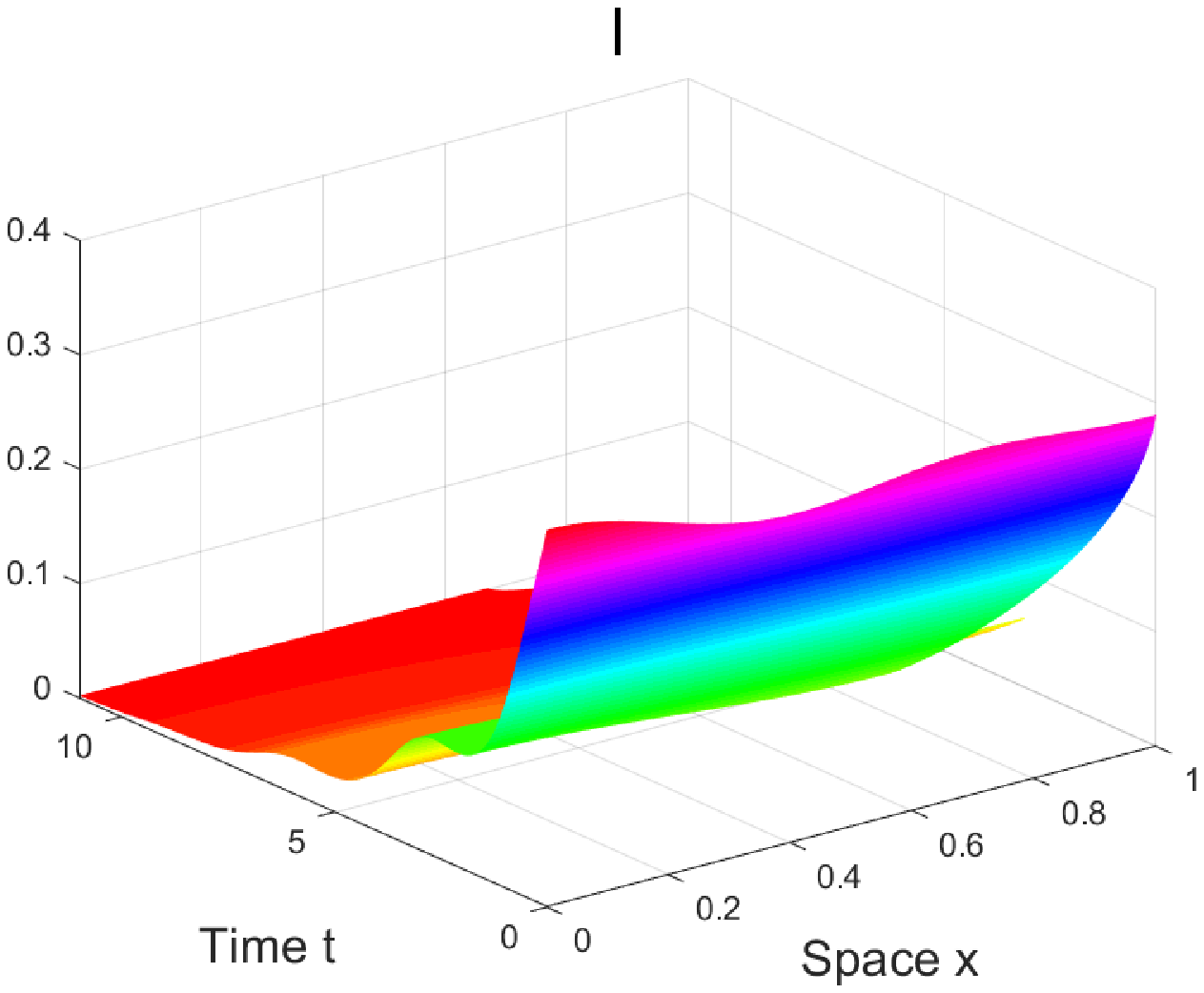}
} }
\subfigure[]{ {
\includegraphics[width=0.28\textwidth]{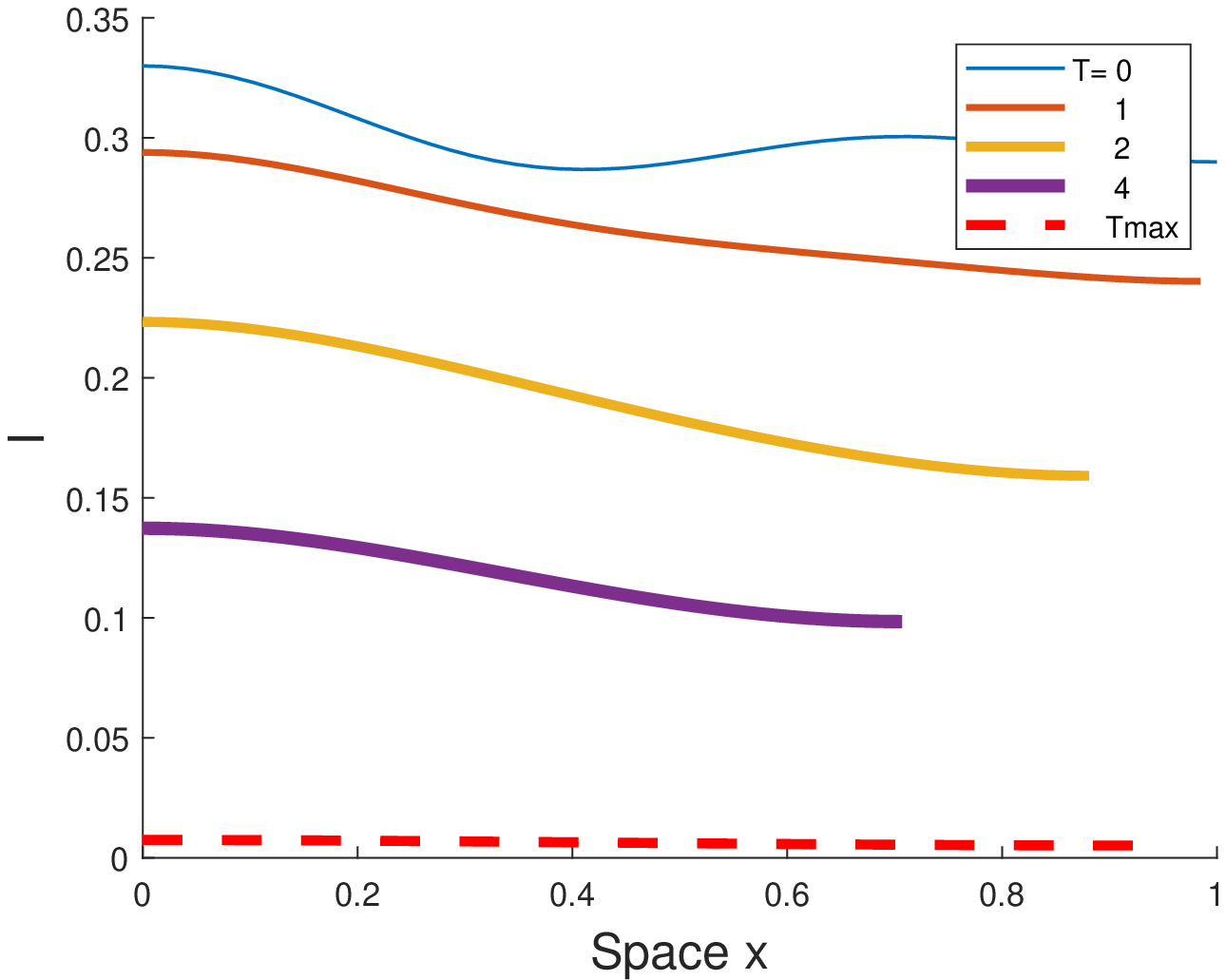}
} }
\subfigure[]{ {
\includegraphics[width=0.28\textwidth]{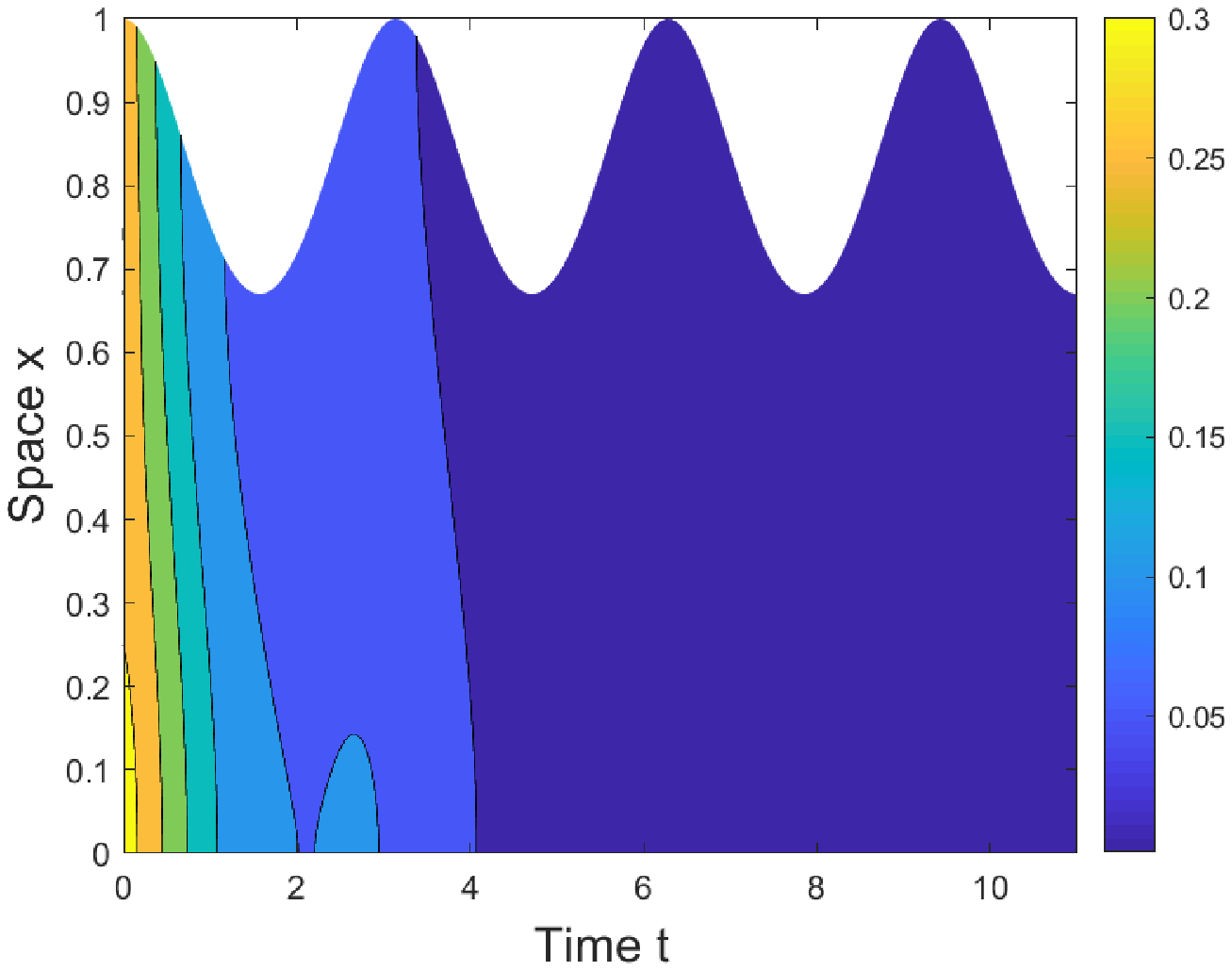}
} }
\caption{\scriptsize $\gamma=0.35+0.88\rho_{6}(t)y$ and $L=1$, ${\mathcal{R}}_0<1$. Graph $(a)$ shows that
infected individual $I$ decays to zero. Graphs $(b)$ and $(c)$, which are the cross-sectional view and
contour map respectively.
}
\label{tu6}
\end{figure}

\begin{figure}[ht]
\centering
\subfigure[]{ {
\includegraphics[width=0.28\textwidth]{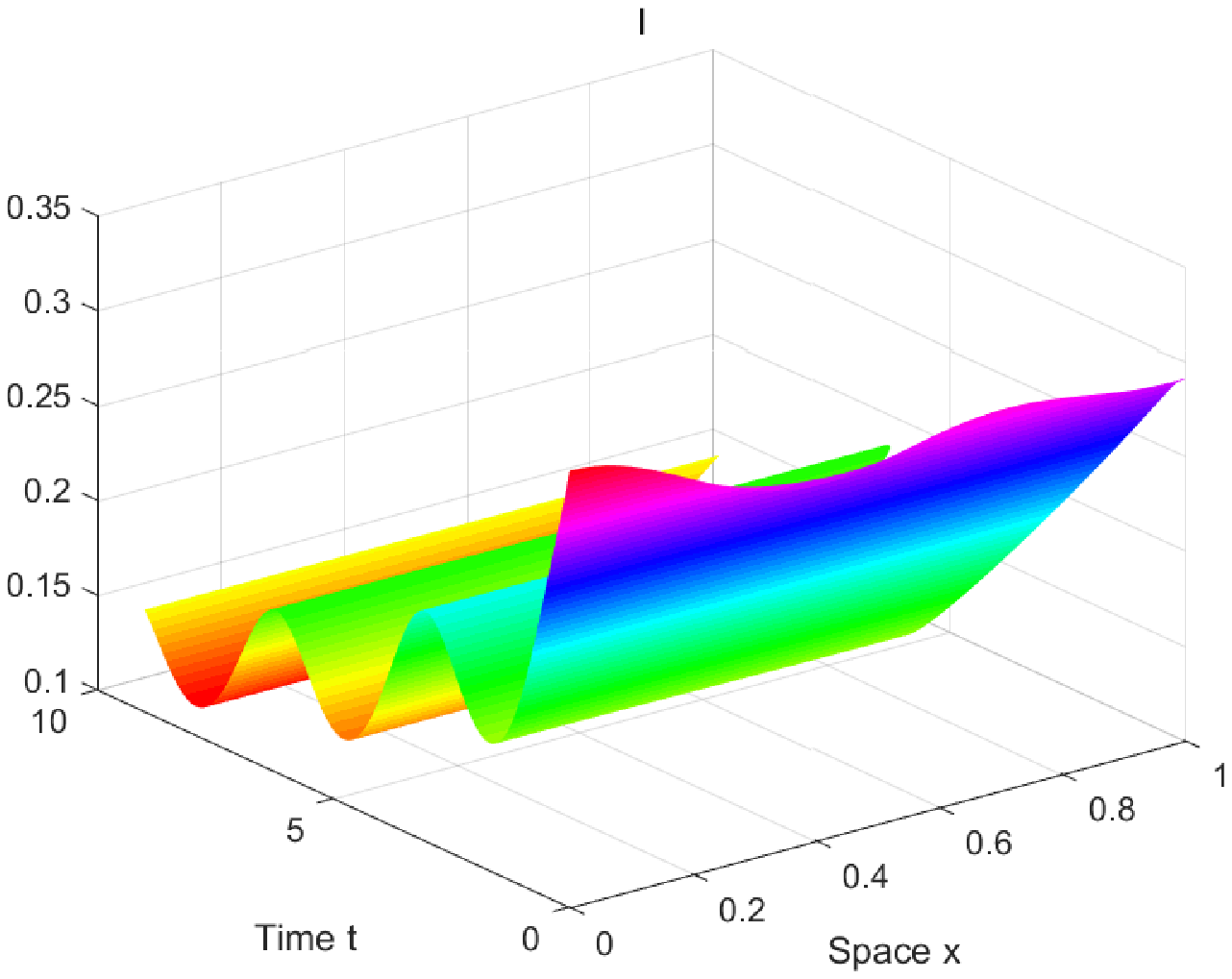}
} }
\subfigure[]{ {
\includegraphics[width=0.28\textwidth]{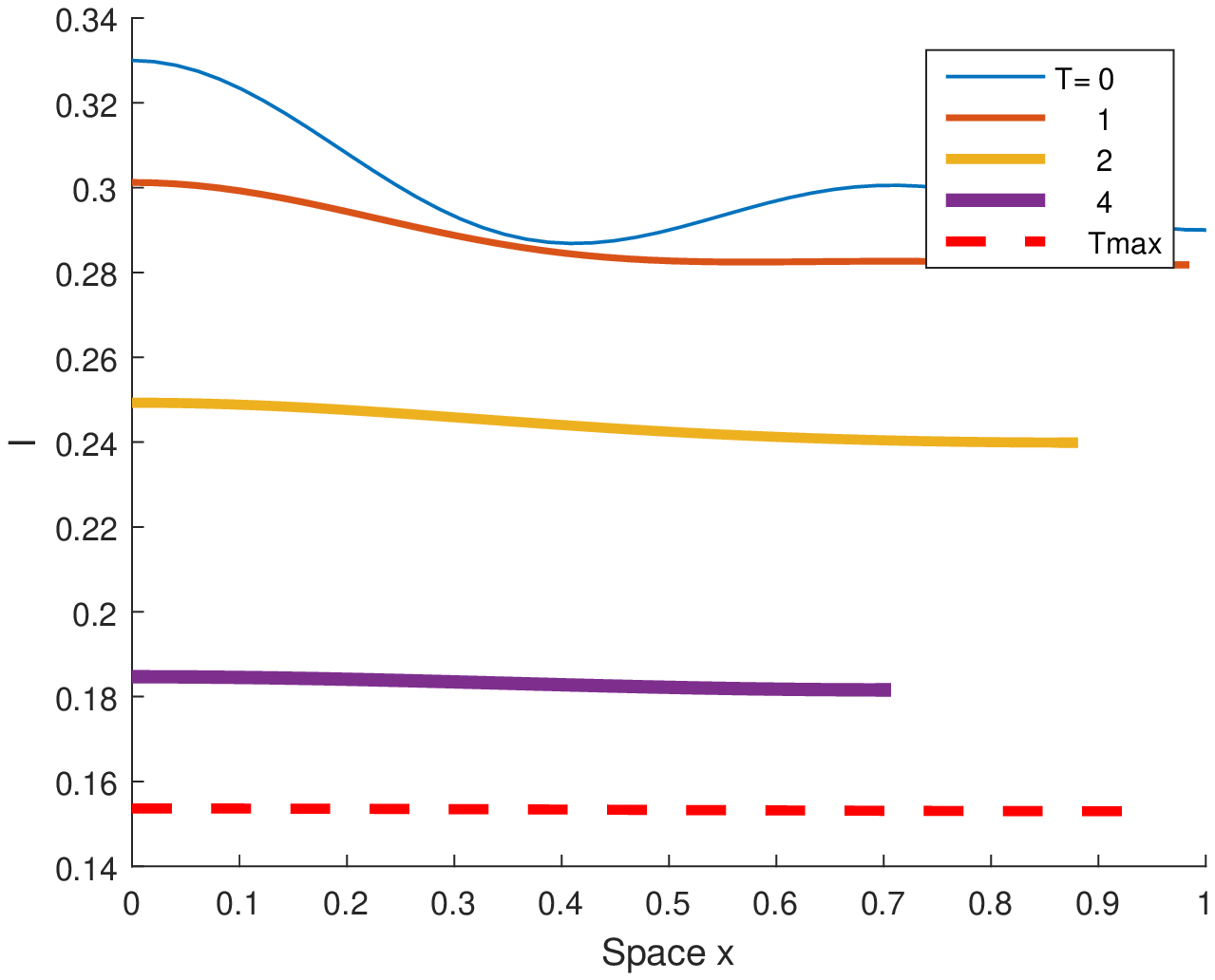}
} }
\subfigure[]{ {
\includegraphics[width=0.28\textwidth]{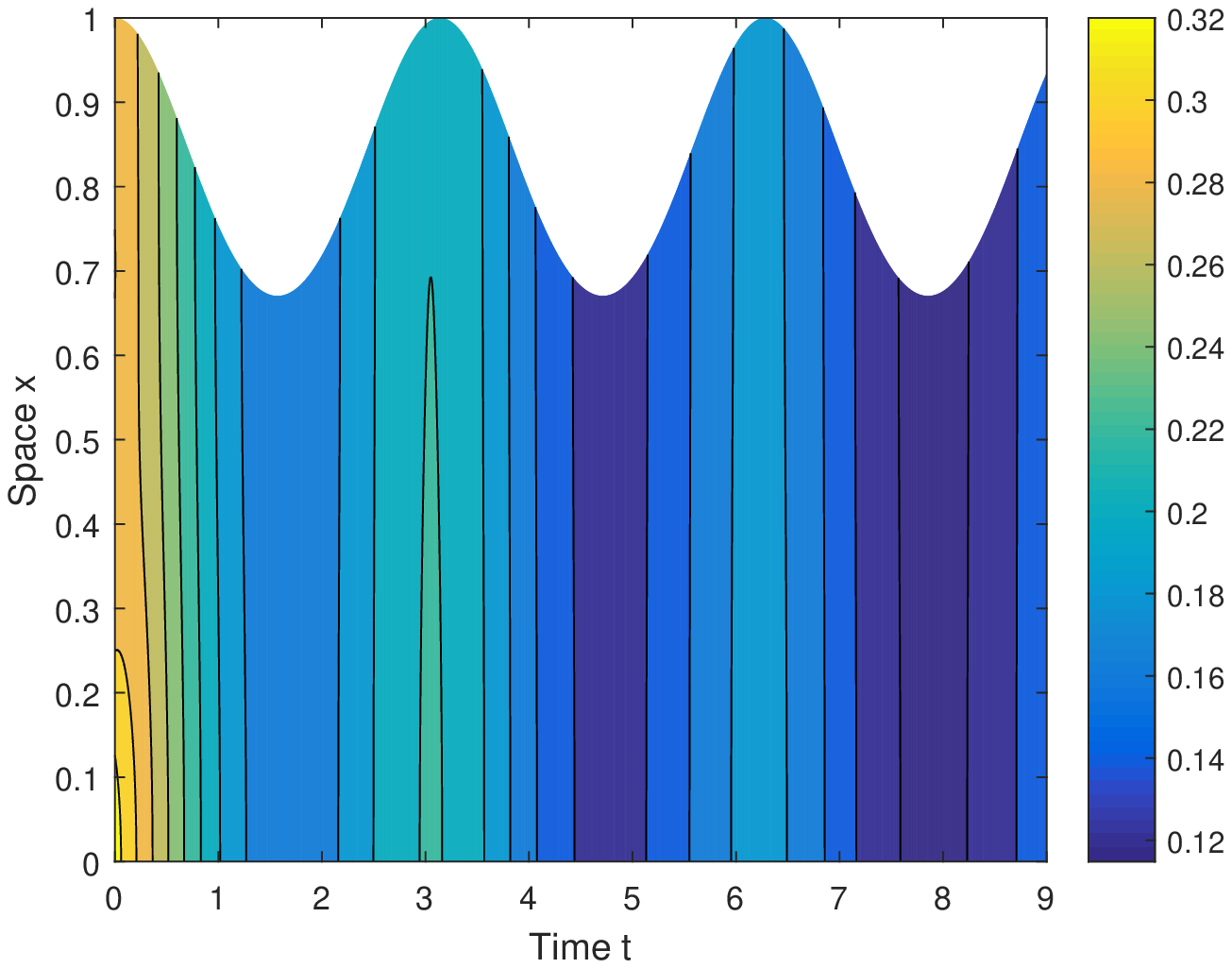}
} }
\caption{\scriptsize $\gamma=0.31+0.01\rho_{6}(t)y$ and $L=2$, ${\mathcal{R}}_0>1$. Graph $(a)$ shows that infected
individual $I$ tends to a positive periodic steady state. Graphs $(b)$ and $(c)$, which are the cross-sectional view
and contour map respectively.
}
\label{tu7}
\end{figure}

\begin{exm}
We first take $\gamma=0.35+0.88\rho_{6}(t)y$ and $L=1$. Due to \eqref{c08}, straightforward calculation alleges

$$
{\mathcal{R}}_0(\rho_6)\leq\frac{\int_0^\pi\max\limits_{y\in\overline{\Omega(0)}}\beta(\rho(t)y)dt}
{\int_0^\pi\min\limits_{y\in\overline{\Omega(0)}}\gamma(\rho(t)y)dt} \approx
\frac{1.0679}{3.3857}<1.
$$
It is easily seen from Fig. \ref{tu6} that $I$ decays to zero.

We now choose $\gamma=0.31+0.01\rho_{6}(t)y$ and $L=2$. It follows from \eqref{c08} that

$$
{\mathcal{R}}_0(\rho_6)\geq \frac{\int_0^\pi\min\limits_{y\in\overline{\Omega(0)}}\beta(\rho(t)y)dt}
{\int_0^\pi\max\limits_{y\in\overline{\Omega(0)}}\gamma(\rho(t)y)dt}\approx  \frac{1.0681}{0.9739}>1.
$$
It is easily seen from Fig. \ref{tu7} that $I$ stabilizes to a positive periodic steady state.
\end{exm}

The evolution of disease-infected domain is an interesting subject and plays an important role in the study of
infectious disease. In order to explore the influence of the periodic evolution in domain on the prevention
and control of the infectious disease, we study the SIS reaction-diffusion model with logistic term on a
periodically evolving domain.

The purpose of this paper is to analyze the effects of diffusion coefficient of infected individuals $d_I$, interval
length $L$ and evolving rate $\rho(t)$ on ${\mathcal{R}}_0$, and ${\mathcal{R}}_0$ as threshold can be used to characterize
stability of the disease-free equilibrium. Firstly, the basic reproduction number ${\mathcal{R}}_0$ is given by
the next generation infection operator and relies on the average value of domain evolution rate $\rho(t)$, see
the formula \eqref{d05}. We then discuss the monotonicity of ${\mathcal{R}}_0$ with respect to diffusion coefficient
$d_I$ and interval length $L$ (see Theorem \ref{main1}), and limiting behavior of ${\mathcal{R}}_0$ if $d_I$ or $L$ is
sufficiently small or sufficiently large in one-dimensional space, see Theorems \ref{main2} and \ref{main3}. We
find that when ${\mathcal{R}}_0<1$, the disease-free equilibrium $(S^{*}(y,t),0)$ is globally asymptotically stable
for system \eqref{a07}-\eqref{a08}, that is to say, the epidemic eventually disappeared; if ${\mathcal{R}}_0>1$,
the disease-free equilibrium $(S^{*}(y,t),0)$ is unstable, in other words, the epidemic persists uniformly, see
Theorem \ref{main}.

In addition, our numerical simulations explain that large evolution rate is unfavorable to control of infectious disease
(see Figs. 1 and 2) and small evolution rate is benefit to control of infectious disease (see Figs. 3 and 4). Meanwhile,
small diffusion coefficient of infected individuals $d_I$ has a positive effect on the presence of infectious
disease (see Figs. 2 and 5). Besides, small habitats facilitate the control of infectious diseases (see Figs. 6 and 7), which
is consistent with our understanding.

\medskip
\medskip

\end{document}